\documentclass[aap,MSNbibl,seceqn,citesort,dvips]{arximspdf}
\usepackage{graphicx}

\doi{10.1214/12-AAP854} 
\volume{23}
\issue{2}
\pubyear{2013}
\firstpage{722}
\lastpage{763}

\makeatletter

\renewcommand{\widetilde}{\tilde}
\renewcommand{\widehat}{\hat}

\newcommand{\varprojlim}{{\displaystyle\mathop{\lim}_{\longleftarrow}}}

\newtheorem{theorem}{Theorem}[section]
\newtheorem{lemma}[theorem]{Lemma}
\newtheorem{corollary}[theorem]{Corollary}
\newtheorem{prop}[theorem]{Proposition}

\newproclaim{remark}[theorem]{Remark}
\newproclaim{example}[theorem]{Example}

\newcommand{\eps}{\varepsilon}
\newcommand{\R}{\mathbb R}

\newcommand{\N}{\mathbb N}
\newcommand{\X}{\mathcal V}
\newcommand{\mmid}{\Vert}
\renewcommand{\mid}{|}

\makeatother

\begin{document}
\begin{frontmatter}

\title{Large deviations for the degree structure in preferential
attachment schemes\thanksref{T1}}
\runtitle{LDP for preferential attachment schemes}

\thankstext{T1}{Supported in part by NSF-09-06713 and NSA Q H982301010180.}

\begin{aug}
\author[A]{\fnms{Jihyeok} \snm{Choi}\ead[label=e1]{jchoi46@syr.edu}}
\and
\author[B]{\fnms{Sunder} \snm{Sethuraman}\corref{}\ead[label=e2]{sethuram@math.arizona.edu}}
\runauthor{J. Choi and S. Sethuraman}
\affiliation{Syracuse University and University of Arizona}
\address[A]{Department of Mathematics\\
Syracuse University\\
Syracuse, New York 13244\\
USA\\
\printead{e1}}
\address[B]{Department of Mathematics\\
University of Arizona\\
617 N. Santa Rita Ave.\\
Tucson, Arizona 85721\\
USA\\
\printead{e2}}
\end{aug}

\received{\smonth{1} \syear{2012}}

%
\begin{abstract}
Preferential attachment schemes, where the selection mechanism is
linear and possibly time-dependent, are considered, and an
infinite-dimensional large deviation principle for the sample path
evolution of the empirical degree distribution is found by
Dupuis--Ellis-type methods. Interestingly, the rate function, which can
be evaluated, contains a term which accounts for the cost of assigning
a fraction of the total degree to an ``infinite'' degree component, that
is, when an atypical ``condensation'' effect occurs with respect to the
degree structure.\looseness=-1

As a consequence of the large deviation results, a sample path a.s. law
of large numbers for the degree distribution is deduced in terms of a
coupled system of ODEs from which power law bounds for the limiting
degree distribution are given. However, by analyzing the rate function,
one can see that the process can deviate to a variety of atypical
nonpower law distributions with finite cost, including distributions
typically associated with sub and superlinear selection models.\looseness=-1
\end{abstract}

%

\end{frontmatter}

\section{Introduction and results}\label{intro}
Preferential attachment processes are graph networks which evolve in
time by linking at each time step a new node to a vertex in the
existing graph with probability based on a selection function of the
vertex's connectivity. Such schemes have a long history in various
guises going back to~\cite{Simon} and~\cite{Yule}; cf. surveys
\cite{Mitzenmacher,Simkin}. More recently, Barab\'{a}si and Albert (BA)
in~\cite{Albert-Barabasi-99} proposed that versions of these processes,
where the selection function is an increasing function of the
connectivity, may serve as models for growing real-world networks such
as the world wide internet web, and types of social structures.

For instance, in a ``friend network,'' a newcomer may have a
predilection to link or
become friends with an individual with high connectivity, or in other
words, one who\vadjust{\goodbreak}
already has many friends. An important property of such reinforcing
networks is
that when the selection function is in a linear form,
asymptotically as time grows, the proportions of nodes with degrees $1,
2, \ldots,
k, \ldots$ converge to a power-law distribution $\langle q(k)\dvtx k\geq1
\rangle$ where $0<\lim_{k\uparrow\infty} q(k)k^{\theta} <\infty$ for
some $\theta>0$. We will say that a network with such a law
of large numbers (LLN) property is ``scale free.'' Since it has been
observed that
the sampled empirical degree structure in many real-world networks has
a ``scale-free'' form, such preferential attachment processes have
become quite popular in
several ways; see \cite
{Albert-Barabasi-02,Barabasi,BBV,BS,Cald,CL06,CH,DM,Du06,Mitzenmacher,Newman,Newman10,NW}
and references therein.

To illustrate more clearly the possible phenomena, consider the
following basic example.\vspace*{-3pt}
%
%
\begin{example}
\label{example1}
Initially, at time $n=1$, the network $G_1$ is composed of two
vertices with a
single (undirected) edge between them. At time $n=2$, a new vertex is
attached to one of the two
vertices in $G_1$ with probability proportional to a function of its
degree to
form the new network $G_2$. This scheme continues: more precisely, at
time $n+1$,
a new node is linked to vertex $x\in G_n$ with probability proportional to
$w(d_x(n))$, that is, chance $w(d_x(n))/\sum_{y\in G_n}w(d_y(n))$,
where $d_z(n)$ is
the degree at time $n$ of vertex $z$, and $w=w(d)\dvtx\N\rightarrow\R_+$
is the
selection function.

In this way, since the initial graph is a tree, all later networks
$G_n$ for $n\geq
1$ are also trees.
Let now $\mathcal{Z}_k(n)$ be the number of vertices in $G_n$ with $k$ links,
$\mathcal{Z}_k(n) = \sum_{y\in G_n} 1(d_y(n) = k)$. We now describe a
trichotomy of growth behaviors corresponding to the
strength and type of
the selection function~$w$~\cite{KR}.

First, when $w$ is linear, say $w(d)=d +\alpha$ for $\alpha>-1$, the
system is ``scale-free.'' As is well
understood in the literature (cf.~\cite{Du06}, Chapter 4), the mean values
$\langle M_k(n)=E[\mathcal{Z}_k(n)]\dvtx k\geq1\rangle$ satisfy rate
equations in the time index
$n\geq1$ which can be solved to show $\lim_{n\uparrow\infty}M_k(n)/n
= q(k)$ for $k\geq1$ where $q$ is in power-law form with $\theta=
3+\alpha$.

Later, in~\cite{BRST01}, when $\alpha=0$, a concentration inequality
was used to show convergence in probability,
$\lim_{n\uparrow\infty}\mathcal{Z}_k(n)/n = q(k)$ for $k\geq1$. We
will call the $\alpha=0$ model the ``classical BA process'' as it was
the model originally analyzed in~\cite{Albert-Barabasi-99}. Also, for
all $\alpha>-1$, P\'olya urn/martingale ideas, and embeddings into
branching processes have given alternative proofs which yield a.s.
convergence; see \mbox{\cite{Mori-01,Rudas-Toth-Valko-07,AGS08}}.

However, in the sublinear case, when $w(d)=d^r$ for $0<r<1$,
although it was shown that a.s. $\lim_{n\uparrow\infty}\mathcal
{Z}_k(n)/n = q(k)$, this
LLN limit $q$ is not a
power law, but in stretched exponential form \cite
{KR,Rudas-Toth-Valko-07}: for $k\geq1$,
%
%
\begin{eqnarray}
\label{stretched}
&\displaystyle q(k) = \frac{\mu}{k^r}\prod_{j=1}^k \biggl(1+\frac{\mu
}{j^r}
\biggr)^{-1}\quad\mbox{and}&\nonumber\\[-9pt]\\[-9pt]
&\displaystyle \mu\mbox{ is determined by } 1 = \sum
_{k=1}^\infty\prod_{j=1}^k
\biggl(1+\frac{\mu}{j^r}\biggr)^{-1}.&\nonumber
\end{eqnarray}
Asymptotically, $\log q(k) \sim-(\mu/(1-r))k^{1-r}$ as
$k\uparrow\infty$. On the other hand, when $r=0$, the case of uniform
attachment when an old vertex is
selected uniformly, an a.s. LLN can also be similarly obtained where
$q$ is geometric. $q(k) = 2^{-k}$ for
$k\geq1$.

In the superlinear case, when $w(d)=d^r$ for $r>1$, ``explosion'' or a
sort of
``condensation effect''' happens in that in the limiting graph a random
single vertex dominates in accumulating connections.
In particular, the
limiting graph is a tree where there is a single random vertex with an infinite
number of children; all other vertices have bounded degree, and of
these only a
finite number have degree strictly larger than $r/(r-1)$; cf., for a
more precise
description,~\cite{Oliveira-Spencer05,KR}. Moreover, a LLN limit,
$\lim_{n\uparrow\infty}E\mathcal{Z}_k(n)/n = q(k)$, is argued where
$q$ is
degenerate in that $q(1)=1$ but $q(k)=0$ for $k\geq2$; cf.~\cite{KR}
and~\cite{Du06}, Chapter 4. Such a limit implies, in the superlinear selection
process, that most of the nodes at step $n$ are leaves.
\end{example}

Since the work of Barab\'{a}si and Albert~\cite{Albert-Barabasi-99},
much effort has been devoted to understand the
degree and other structures in generalized versions of these graphs.
A partial selection of this large literature includes: more on degree structure
\cite{DMpaper,DKM,Fort,FVC,Jans,KRR,KR-finite}; growth and location of
the maximum degree
\cite{Mori-05,AGS08,dereich-2009}; spectral gap and cover time of a
random walk on
the graph~\cite{MPS,Cooper-Frieze07}; width and diameter
\cite{Katona,Bollobas-Riordan,Dereich-Monch-Morters}; graph limits
\cite{Berger-Borgs-Chayes-Saberi-09,Bhamidi,BCLSV,RS}.\vspace*{8pt}

\textit{Connection between urns and degree structure.}
If, however, one focuses only on the degree structure of the growing network,
then it may be helpful to view the degree distribution evolution in
terms of
``balls-in-bins'' or ``P\'olya urn'' models. For instance, in the
previous example, every
new connection that a vertex gains can be represented by a new ball
added to a
corresponding urn in a collection of urns.
More precisely, at time $n=1$, there are two urns, each possessing one
ball, in the initial collection
$U(0)$. At time $j+1$, a new urn with one ball is included in the
collection, and also one
ball is added to an existing urn $x\in U(j)$ with probability
proportional to $w(b_x)$ where $b_x$ is the number of balls in urn $x$.
Then, $\mathcal{Z}_k(n)$
translates to the number of urns in $U(n)$ with $k$ balls for $k\geq1$.

A comprehensive form of such an urn model was formulated by Chung,
Handjani and Jungreis (CHJ) in~\cite{Chung-Handjani-Jungreis},
motivated by the work in
\cite{DEM} and~\cite{KR} on the
organization of web tree-graph models. See also~\cite{BBCS-05}
and~\cite{Mori-05} for other work connecting urns to degree structure.

The CHJ model is as follows. Given an initial finite collection
of urns each containing one ball, at subsequent times, with probability
$p$, a new
urn with one ball is created and added to the collection or, with
probability $1-p$,
a~new ball is put in one of the existing urns $x$ with probability\vadjust{\goodbreak}
proportional to
$(b_x)^r$ where $b_x$ is the number of balls in $x$.
It was
proved in~\cite{Chung-Handjani-Jungreis}, among other results, when
$r=1$ and $p>0$, analogous to linear selection preferential attachment
graphs, that the
empirical urn size distribution converges to a power law with $\theta
= 1 + (1-p)^{-1}$.\vspace*{8pt}

In this context, our purpose is to study a generalized preferential
attachment process of urns, where at each time step a new urn is created
and a new ball is added to it or an existing urn according to a time dependent
linear selection function, which includes the evolving degree
structure of linear selection
preferential attachment model discussed above, and also
a version of
the $r=1$ CHJ urn model~\cite{Chung-Handjani-Jungreis}. We defer to
Section~\ref{model} the exact description of our scheme.

As mentioned in~\cite{Drinea-Frieze-Mitzenmacher},
understanding preferential attachment or urn systems, where the
selection function
depends on time, allows for more realistic models
given real world networks are time-dependent.
However, it appears most of the work on time-dependent schemes
consists of rate equation formulations
(\cite{Dorogovtsev-Mendes}, Section E),~\cite{KR} and related work, in
models where at
each step a random number of links or balls may be added to the
structure~\cite{Cooper-Frieze-general,AGS08}.

Given this background, detailing the large deviation behavior of the degree
distribution in time-dependent preferential attachment schemes is a
natural problem
which gives much understanding of typical and, in particular, atypical
evolutions. We
remark, even in the usual time-homogeneous models, large deviations of
the degree
structure is an open question.

Previous large deviation work in preferential attachment models has
focused on one-dimensional objects, for instance, the number of leaves
\cite{Bryc-Minda-S}, or the degree growth of a single vertex with
respect to dynamics where any vertex may attach to a newly added vertex
with a small chance~\cite{dereich-2009}. See references cited in
\cite{Bryc-Minda-S} for large deviations work with respect to other
types of random trees and balls-in-bins models.

Our main work in this article includes an infinite-dimensional sample
path large deviation principle
(LDP) for an array of empirical urn ball size distributions $\{
\langle Z^n_k(j)/n\dvtx k\geq0\rangle\dvtx0\leq j\leq
n\}_{n\geq1}$, when the initial configuration, not necessarily fixed,
satisfies a
limit condition (Theorem~\ref{LDPinfin}). Here, $Z^n_k(j)$ stands for
the count of urns with $k$
balls at time $j$ in the $n$th row of the array.
Part of these results is a finite-dimensional LDP with respect to the
numbers of
urns with less than $d$ balls for $d<\infty$ (Theorem~\ref{LDPmain}).

As an application of the large deviations results, we obtain an a.s.
sample path
LLN for the urn counts in terms of a system of coupled ODEs (Corollary
\ref{LLN}), which, for homogeneous schemes complements fixed time LLNs
mentioned
earlier, and gives a different way to derive them aside from the rate equation
method mentioned in Example~\ref{example1}. Finally, the LLN limit trajectories
are shown to have power law-type behavior in terms of bounds (Corollary
\ref{powerlaw}), although the general behavior can interpolate between these
bounds; see
Figure~\ref{fig1}.\vadjust{\goodbreak}

Interestingly, the infinite-dimensional rate function $I^\infty$ can be
calculated
on scaled urn ball size path distributions $\xi= \{\langle\xi_k(t)\dvtx
k\geq0\rangle\dvtx0\leq t\leq1\}$.
Here, since in our model, exactly one ball is added to the urn
collection at
each microscopic time, $\xi_k(t)/(t+c)$ is the fraction of urns with
size $k$
at macroscopic time $t\geq0$ where $c= \sum_{k\geq0}\xi_k(0)$ is the
initial mass, that is
the scaled initial number of urns.
It is natural then to ask which trajectories $\xi$ have finite cost,
$I^\infty(\xi)<\infty$.

It turns out ``no mass can be lost,'' that is, all finite cost
paths $\xi$ are such that the proportions $\{\xi_k(t)/(t+c)\}_{k\geq
0}$ form a
probability distribution, $\sum_{k\geq0} \xi_k(t)/(t+c) \equiv1$.
Also, a variety
of nonpower law distributions can be achieved with finite rate at any
time $0<t\leq1$, including
geometric and stretched exponential distributions discussed in Example
\ref{example1}.

Intriguingly, on the other hand, ``some of the weight may be lost'' in
certain finite
rate trajectories, that is, the scaled mean urn ball size of the system
may satisfy
a ``weight loss'' property at a time $0<t\leq1$, $\sum_{k\geq0} k\xi
_k(t)/(t+\tilde{c}) <1$, where $\tilde{c} =
\sum_{k\geq0} k\xi_k(0)$ is
the scaled initial total urn ball size, even though the pre-limit
quantity equals $1$ at
all steps in the urn growth scheme.
We dub a trajectory $\xi$ with this ``weight
loss'' property at some time $0<t\leq1$ as being ``condensed.'' For
instance, a ``condensed'' path arises
when $c=\tilde c=0$ and a finite number of the urns take in eventually
all the
balls. In this case, almost all the urns created are empty, and the
associated path
satisfies $\xi_0(t) = t$ for $0\leq t\leq1$, $\xi_k(t) \equiv0$ for
$k\geq1$, and $\sum_kk\xi_k(t) \equiv
0$. It turns out this path, associated with superlinear selection
preferential attachment models (cf. Example~\ref{example1}), has
finite cost.

Moreover, the rate function $I^\infty$ contains a term which measures
the cost of
``condensation'' when some of the flow of urn ball size in the scaling
limit escapes
toward urns with ``infinite'' size. In addition, we point out, at any
time $0\leq t\leq1$,
LLN distributions arising
from either sublinear or superlinear selection preferential attachment
models may be achieved with finite cost.
One might interpret that although the
linear selection process is typically ``scale-free,'' since it is
between, in a sense, sublinear and
superlinear selection models, its atypical degree distribution
structure may
include the typical behavior of its sub and superlinear relatives.
See Remark~\ref{remarkinfinity} and Example~\ref{example2} for more
details and discussion.

We also remark that the large deviations and other work are, with
respect to the process,
starting from either ``small'' or ``large'' initial configurations, that
is, when the
initial urn collection has $o(n)$ balls (e.g., finite), or when
the size of
the collection is on order $n$, respectively. It appears these initial
configurations, which enter into all result statements, have not been considered
before, in general.

The main idea for the results is to extend a variational control problem/weak
convergence approach of Dupuis and Ellis (cf.~\cite{Dupuis-Ellis}) to establish
finite-dimensional LDPs in the time-dependent setting. Then, a
projective limit
approach, and some analysis to identify the rate function, is used to
obtain the
infinite-dimensional LDP. For the LLN and power-law corollaries, a
coupled system
of ODEs, which governs the typical degree distribution evolution, is
identified,
and analyzed.

To be concrete, we have focused upon models where the network is
incremented by one urn and one ball each time, which include basic
models. However, the methods here should be of use to analyze the large
deviations of the degree structure in other combinatorial models with
different increment structure: for instance, the evolving graph model
discussed in~\cite{CL06}, Chapter~3, where at each time with
probability $p$ a new vertex is preferentially attached to an old one,
and with probability $1-p$, an edge is added between two old nodes
selected preferentially, and the BA graphs where, instead of only one
vertex, $m\geq2$ vertices are introduced and preferentially connected
at each time; cf.~\cite{Du06}, Chapter 4.

\subsection{Model}\label{model}

Let $p(t)\dvtx[0,1]\to[0,1]$ and $\beta(t)\dvtx[0,1]\to[0,\infty)$ be given
functions. We define an urn configuration $U$ as a finite collection of
urns, each urn $x\in U$ containing a nonnegative number of balls $b_x$.
We now specify an evolving array $\{U^n(j)\dvtx0\leq j\leq n\}_{n\geq1}$
of urn configurations by the following time-dependent iterative scheme.
In the $n$th row of the array:
\begin{itemize}
\item Start at step $0$, with a given initial urn configuration $U^n(0)$.
\item At step $j+1\leq n$, to form a new urn configuration $U^n(j+1)$,
we first
create and include a new urn with no ball. Then:

\begin{itemize}
\item with probability $p(j/n)$, we place a new ball in this urn;
\item with probability $1-p(j/n)$, we place a new ball in one of the
other urns
$x\in U^n(j)$ with probability
\[
\frac{b_x+\beta(j/n)}{\sum_{y\in U^n(j)} (b_y+\beta(j/n))}.
\]
\end{itemize}

\end{itemize}
We will call, for urn $x\in U^n(j)$, the term $b_x + \beta(j/n)$ as the
``weight'' of the urn
at time $j$ in the $n$th row of the process.
Let now $|U^n(j)|$ and $B^n(j) = \sum_{x\in U^n(j)} b_x$ be the total
number of urns and balls in $U^n(j)$, respectively.
Then, the number of urns $|U^n(j)| =
|U^n(0)| + j$ and the total number of
balls $B^n(j)=
B^n(0) +j$. Also, the total weight of the configuration at time $j$ is
\[
s^n(j):= \sum_{y\in U^n(j)} \bigl(b_y + \beta(j/n) \bigr)
=
\bigl(1+\beta
(j/n)\bigr)j +
B^n(0) + \beta(j/n)|U^n(0)|.
\]

The above urn scheme, as discussed in the \hyperref[intro]{Introduction}, may be viewed in
terms of the evolving degree structure in a preferential attachment
random graph
process with time-dependent selection function $w(d;j,n) = d + \beta
(j/n)$. Here,
the step of including a new empty urn and incrementing the number of
balls in an old urn\vadjust{\goodbreak}
corresponds to an edge being placed between a new node, with degree
$1$, and an old
vertex in the existing graph whose degree is consequently incremented. In
particular, when $p$ and $\beta$ are in particular forms, we recover
the following
models:

\begin{longlist}[(3)]
\item[(1)] ``Classical'' BA process. When $p(t)\equiv0$ and
$\beta(t)\equiv1$, the scheme is time-homogeneous. When the initial urn
configuration consists of two empty urns, the probability of selecting
an urn $x$
with $k\geq0$ balls at time $j\geq0$ is $(k+1)/(2(j+1))$, which
matches the
selection process in the evolution of the degree structure in the BA
preferential
attachment graph scheme at times $j+1\geq1$ with selection function
$w(d) = d$, as
discussed in Example~\ref{example1}, where urns with $k\geq0$ balls
correspond to
vertices with degree $d=k+1\geq1$.

\item[(2)] ``Offset'' BA processes. When $p(t)\equiv0$ and
$\beta(t)\equiv\beta\geq0$, again the scheme is time-homogeneous, and
urns with
$k\geq0$ balls correspond to vertices with degree $k+1\geq1$.
However, now the
weight of an urn with $k$ balls is $k+\beta$, in a sense ``offset'' from
the classical
BA scheme. Correspondingly, the urn selection scheme is the same as the growth
process of the degree structure in the preferential attachment model
with selection
function $w(d)=d+\alpha$ with $\alpha= \beta-1$ as specified in Example
\ref{example1}.

\item[(3)] CHJ model of P\'olya urns. When
$p(t)\equiv p$ and $\beta(t)\equiv\beta\geq0$,
the evolution of the number of urns of size
$k\geq0$ corresponds to
a version of the $r=1$ CHJ model discussed in the \hyperref[intro]{Introduction}.
However, we
note, in our
model, an empty urn is added at each step with probability $1-p$, and
these empty urns are kept track of in our scheme.
When $\beta=0$, the dynamics of urns of size $k\geq1$ is the $r=1$
CHJ model since the empty
urns have no weight, and once created, they cannot be selected to fill
in later steps,
and do not influence the structure of urns with $k\geq1$ balls.
\end{longlist}

For $n\geq1$, let ${Z^n_i(j)}$ be the number
of urns in the $n$th row of the urn array process with $i\geq0$ balls
at time $0\leq j\leq n$ and, for $d\geq0$, let $\bar
Z_{d+1}^n(j)$
denote the number of urns with more than $d$ balls at time $0\leq j\leq n$.
These quantities satisfy
\begin{eqnarray*}
\sum_{i=0}^{d} Z_{i}^n(j) + \bar Z_{d+1}^n(j) & =& |U^n(0)|+j,\\
\sum_{i=0}^{d} i Z_{i}^n(j) + (d+1)\bar Z_{d+1}^n(j) & \leq& B^n(0) +j.
\end{eqnarray*}

Define now vectors in $\R^{d+2}$,
%
\begin{eqnarray}
\mathbf f^d_0&:=& \langle0,1,0,\ldots,0\rangle,\qquad
\mathbf f^d_i:= \langle1,0,\ldots,0,-1,1,0,\ldots,0\rangle\nonumber\\
&&\eqntext{\mbox{where } {-1} \mbox{ is at the
}(i+1)\mbox{th position for } 1\leq i \leq d,}\\
\mathbf f^d_{d+1}&:=& \langle1,0,\ldots,0\rangle.\nonumber
\end{eqnarray}
For $\mathbf y=\langle y_0,\ldots, y_d, y_{d+1}\rangle\in\R^{d+2}$ and
$0\leq
i\leq d+1$, denote
\[
[\mathbf y]_i:=\sum_{l=0}^i y_l.
\]
Note that
%
%
\begin{eqnarray}\label{increments}
&\displaystyle 0\leq[\mathbf{f}^d]_i\leq1 \qquad\mbox{for $0\leq i\leq d$},&
\nonumber\\[-8pt]\\[-8pt]
&\displaystyle [\mathbf{f}^d]_{d+1}=
1\quad \mbox{and}\quad 0\leq\sum_{i=0}^{d+1}
(1-[\mathbf{f}^d]_{i})\leq1.&\nonumber
\end{eqnarray}

Consider now the ``truncated'' degree distribution
\[
\{\mathbf Z^{n,d}(j):=\langle Z_{0}^n(j),\ldots, Z_{d}^n(j), \bar
Z_{d+1}^n(j)\rangle\mid0\leq j\leq n\},
\]
where $\bar Z_{d+1}^n(j) = \sum_{k\geq d+1}Z_k^n(j) = j+ |U^n(0)| -
\sum_{k=0}^d
Z_k^n(j)$, which forms a discrete time Markov chain with initial state
$\mathbf
Z^{n,d}(0)$ corresponding to the initial urn configuration $U^n(0)$ and one-step
transition property,
\begin{eqnarray*}
&&\mathbf Z^{n,d}(j+1) - \mathbf Z^{n,d}(j)\\
&&\qquad =
\cases{
\mathbf f^d_0, &\quad with prob. $\displaystyle
p(j/n) + \bigl(1-p(j/n)\bigr)\frac{\beta(j/n)
Z_{0}^{n}(j)}{s^n(j)}$\vspace*{2pt}\cr
&\qquad for $i=0$,\vspace*{2pt}\cr
\mathbf f^d_i, &\quad with prob. $\displaystyle
\bigl(1-p(j/n)\bigr)\frac{(i+\beta(j/n))
Z_{i}^{n}(j)}{s^n(j)}$\vspace*{2pt}\cr
&\qquad for $1\leq i\leq d$,\vspace*{2pt}\cr
\mathbf f^d_{d+1}, &\quad with prob.
$\displaystyle \bigl(1-p(j/n)\bigr)\biggl(1-\frac
{\sum_{i=0}^d
(i+\beta(j/n)) Z_{i}^{n}(j)}{s^n(j)}\biggr)$.}
\end{eqnarray*}
We also define the ``full'' degree distribution
\[
\{\mathbf Z^{n,\infty}(j):=\langle Z_{0}^n(j),\ldots, Z_{d}^n(j),
\ldots\rangle\mid0\leq j\leq n\},
\]
which is also a Markov chain on $\R^\infty$ with increments
\begin{eqnarray*}
&&\mathbf Z^{n,\infty}(j+1) - \mathbf Z^{n,\infty}(j)\\
&&\qquad=
\cases{
\mathbf f^\infty_0, &\quad with prob.
$\displaystyle p(j/n) + \bigl(1-p(j/n)\bigr)\frac
{\beta(j/n)
Z_{0}^{n}(j)}{s^n(j)}$\vspace*{2pt}\cr
&\qquad for $i=0$,\vspace*{2pt}\cr
\mathbf f^\infty_i, &\quad with prob.
$\displaystyle \bigl(1-p(j/n)\bigr)\frac{(i+\beta(j/n))
Z_{i}^{n}(j)}{s^n(j)}$\vspace*{2pt}\cr
&\qquad for $i\geq1$,}
\end{eqnarray*}
where $f_0^\infty= \langle0,1,0,\ldots, 0, \ldots\rangle$ and
$f_i^\infty= \langle
1,0,\ldots,
0,-1,1,0,\ldots, 0, \ldots\rangle$ with the ``$-1$'' being in the
$(i+1)$th place.

We will assume throughout the following initial condition, which
ensures a LLN at time $t=0$.
With respect to constants $c^n_i, c^n, \tilde c^n\geq0$, for $i\geq
0$, define
\[
c_i^n:= \frac1n Z^n_i(0),\qquad c^n:= \sum_{i\geq0}c^n_i
\]
and
\[
\tilde c^n:= \sum_{i\geq0} ic_i^n.
\]

\begin{longlist}[(LIM)]
\item[(LIM)]
For constants $c_i, c, \tilde c\geq0$, we have
\[
c_i:= \lim_{n\uparrow\infty} c_i^n\quad \mbox{and}\quad\tilde c: =
\lim_{n\uparrow\infty}
\tilde c^n = \sum_{i\geq0}ic_i < \infty.
\]
Consequently,
$c:=\lim_{n\uparrow\infty} c^n
= \sum_{i\geq0}c_i <\infty$.
\end{longlist}
In the previous sentence, the $c^n$ limit follows from the
uniform bound, $\sum_{i\geq A} c^n_i \leq A^{-1}\sum_{i\geq0} ic^n_i
\rightarrow
\tilde c/A$.
Define also
\[
\bar c^d:= \sum_{i\geq d+1} c_{i}
\quad\mbox{and}\quad\mathbf c^d
:=\langle c_0,\ldots,
c_d, \bar c^d\rangle.
\]

We remark one can classify the initial configurations depending on when
$c_i\equiv
0$ or when $c_i>0$ for some $i\geq0$.
\begin{itemize}
\item(Small configuration) $c_i\equiv0$ for any $i\geq0$. Here, the
initial urn
configurations are small in that their size is $o(n)$. This is the case
when the
initial configurations do not depend on $n$, for instance.

\item(Large configuration) $c_i>0$ for some $i\geq0$. Here, the
initial state
is already a partly-developed configuration whose size is of order $n$.
\end{itemize}

We also note, when the initial urn configurations correspond to initial tree
configurations in the corresponding preferential attachment process, some
restrictions in the values of $c_i$ arise. One may verify that a
graph with $n$ vertices with degrees $d_1, \ldots, d_n$ is a tree
exactly when $\sum_{i=1}^n d_i = 2(n-1)$. Hence, since in the initial
graph of the $n$th
row, the number of
vertices equals $n\sum_{k\geq0} c^n_k$, and the sum of degrees equals
$n\sum_{k\geq0}(k+1)c^n_k$ (recall the correspondence between urn sizes
and degrees discussed in the \hyperref[intro]{Introduction}), we have $n\sum_{k\geq0}
(k+1)c^n_k = 2(n\sum_{k\geq0}c^n_k -1)$. By (LIM), we have then
$\tilde c
= c$.

In addition, we note (LIM) specifies an initial limiting degree
distribution which has full ``weight'' or in other words is not
``condensed,'' that is, $\tilde c = \lim_{n\uparrow\infty}\tilde
c^n = \sum_{i\geq0} ic_i$. See Remark~\ref{LLNrmk}, however, for
comments when the
initial distribution is ``condensed,'' that is,
$\tilde c > \sum_{i\geq0} ic_i$.

Our results will be on the family of processes
$\mathbf{X}^{n,d}=\{\mathbf{X}^{n,d}(t)\dvtx0\leq t\leq1\}$ and
$\mathbf{X}^{n,\infty}=\{\mathbf{X}^{n,\infty}(t)\dvtx0 \leq t\leq1\}$
obtained by linear interpolation of the discrete-time Markov chains
$\frac
1n\mathbf{Z}^{n,d}(j)$ and $\frac1n\mathbf{Z}^{n,\infty}(j)$, respectively.
For $0\leq t\leq1$, let
\begin{eqnarray*}
\mathbf{X}^{n,d}(t)&:=& \frac1 n \mathbf Z^{n,d}({\lfloor nt\rfloor
})
+
\frac{nt-\lfloor nt\rfloor}{n}\bigl(\mathbf Z^{n,d}({\lfloor
nt\rfloor+1})
- \mathbf Z^{n,d}({\lfloor nt\rfloor})\bigr),\\
\mathbf{X}^{n,\infty}(t)&:=& \frac1 n \mathbf Z^{n,\infty
}({\lfloor
nt\rfloor})
+
\frac{nt-\lfloor nt\rfloor}{n}\bigl(\mathbf Z^{n,\infty}({\lfloor
nt\rfloor+1})
- \mathbf Z^{n,\infty}({\lfloor nt\rfloor})\bigr).
\end{eqnarray*}
The\vspace*{1pt} trajectories $\mathbf X^{n,d}$ lie in
$C([0,1];\R^{d+2})$, and are Lipschitz with constant at most $1$,
satisfying $\mathbf X^{n,d}(0)=\frac 1n\mathbf Z^{n,d}(0)$. On
the\vspace*{1pt} other hand, the infinite distribution $\mathbf
X^{n,\infty}\in \prod_{i=0}^\infty C([0,1]; \R)$, considered with the
product topology, where $\mathbf X^{n,\infty}(0)= \frac1n \mathbf
Z^{n,\infty}(0)$. In both cases, although $\mathbf X^{n,d}(t)$ and
$\mathbf X^{n,\infty}(t)$ are not necessarily probabilities because it
is possible that we do not normalize by the total mass; they are,
however, finite distributions.

We now specify the assumptions on $p(t)$ and $\beta(t)$ used for the
main results.
\begin{longlist}[(ND)]
\item[(ND)]
$p$ and $\beta$ are piecewise continuous and, for constants $p_0$,
$\beta_0$ and $\beta_1$,
\[
0\leq p(\cdot)\leq p_0<1 \quad\mbox{and}\quad 0<\beta_0\leq\beta
(\cdot)<
\beta_1<\infty.
\]
\end{longlist}
We discuss (ND) more in the remark after Theorem~\ref{LDPmain}.

We note, throughout the article, that we use conventions
%
%
\begin{eqnarray}\label{convention}
0\log0 &=&
0/0 = 0\cdot\pm\infty= 1/\infty=0,\nonumber\\
\pm1/0 &=& \pm\infty\quad \mbox{and}\\
E[X;\mathbb A] &=&\int_{\mathbb A} X \,dP.\nonumber
\end{eqnarray}

\subsection{Results}
We now recall the statement of a large deviation principle (LDP).
A sequence $\{X^n\}$ of random variables taking values in a complete separable
metric space $\X$
satisfies the LDP with rate $n$ and good rate function $J\dvtx\X\to
[0,\infty]$
if for each $M<\infty$, the level set $\{x\in\X\mid J(x)\leq M\}$ is
a compact
subset of $\X$,
that is, $J$~has compact level sets, and if the following two
conditions hold:
\begin{longlist}
\item{Large deviation upper bound.} For each closed subset $F$ of $\X$,
\[
\limsup_{n\to\infty} \frac1 n \log P\{X^n\in F\} \leq-\inf
_{x\in
F} J(x).
\]

\item{Large deviation lower bound.} For each open subset $G$ of $\X$,
\[
\liminf_{n\to\infty} \frac1 n \log P\{X^n\in G\} \geq-\inf
_{x\in
G} J(x).
\]
\end{longlist}

For $d\geq0$, we now state the LDP for $\{\mathbf{X}^{n,d}(t)\mid
0\leq t\leq1\}$.
Define the function $I_d\dvtx C([0,1];\R^{d+2})\to[0,\infty]$ given
by
\begin{eqnarray*}
I_d(\varphi)&=& \int_0^1 \bigl(1-[\dot{\varphi}(t)]_0\bigr) \log
\frac{1-[\dot{\varphi}(t)]_0}{p(t)+(1-p(t))\frac{\beta(t)\varphi_0(t)}{(1+\beta(t))t+\tilde
c+c\beta(t)}}
\\
&&\hspace*{12.5pt}{} + \sum_{i=1}^d \bigl(1-[\dot{\varphi}(t)]_i\bigr) \log
\frac{1-[\dot{\varphi}(t)]_i}{(1-p(t))\frac{(i+\beta(t))\varphi_i(t)}{(1+\beta(t))t+\tilde
c+c\beta(t)}}\\
&&\hspace*{12.5pt}{} + \Biggl(1-\sum_{i=0}^d\bigl(1- [\dot{\varphi}(t)]_i\bigr)\Biggr) \log
\frac{1-\sum_{i=0}^d(1-[\dot{\varphi}(t)]_i)}{(1-p(t))(1-\frac{\sum_{i=0}^d
(i+\beta(t))\varphi_i(t)}{(1+\beta(t))t+\tilde c+c\beta(t)})}\,dt,
\end{eqnarray*}
where $\varphi(0)=\mathbf c^d$, $\varphi_i\geq0$ is Lipschitz with
constant 1 such
that $0\leq[{\dot\varphi}(t)]_i\leq1$ for $0\leq i\leq d$, $\sum
_{i=0}^{d+1}
\dot\varphi_i(t)=1$, $\sum_{i=0}^d
(1-[\dot\varphi(t)]_i)=\sum_{i=0}^{d+1}i\dot\varphi_i(t)\leq1$
for almost all $t$,
and the integral converges; otherwise, $I_d(\varphi)=\infty$. It will
turn out that
$I_d$ is convex and is a good rate function.

To explain the last condition in the definition of $I_d$, note that
$\varphi_{d+1}(t)$ represents the fraction of urns with size at least
$d+1$, so that
$(d+1)\varphi_{d+1}(t)$ is the truncated fraction of balls in these
urns. Since the process increments by one ball at each step, it makes
sense to specify $\sum_{i=0}^{d}(1-[\dot\varphi(t)]_i)=\sum
_{i=0}^{d+1}i\dot\varphi_i(t) \leq
1$ or that $\sum_{i=0}^{d+1}i\varphi_i(t) \leq t+ \tilde c$ if
$I_d(\varphi)<\infty$.

The rate function can be understood as follows: in order for $\mathbf
X^{n,d}$ to
deviate to $\varphi$, at time $t$, the process should behave as if the
increment
probabilities $v_i$ of $\mathbf f^d_i$ are such that the mean $\sum_{i=0}^d
v_i\mathbf f^d_i + v_{d+1} \mathbf f^d_{d+1} = \dot\varphi$.
In the proof of Theorem~\ref{LDPmain}, we show $v_i = 1-[\dot\varphi
]_i$ for
$0\leq i\leq d$ and $v_{d+1}=1-\sum_{j=0}^d (1-[\dot\varphi]_j)$.
But, the natural evolution increment probabilities $u_i$, given the
process is in
state $\varphi(t)$, are
$u_0=p(t)+(1-p(t))\frac{\beta(t)\varphi_0(t)}{(1+\beta(t))t+\tilde
c+c\beta(t)}$,
$u_i = (1-p(t))\frac{(i+\beta(t))\varphi_i(t)}{(1+\beta(t))t+\tilde
c+c\beta(t)}$
for $1\leq i\leq d$ and $u_{d+1}=(1-p(t))(1-\frac{\sum_{i=0}^d
(i+\beta(t))\varphi_i(t)}{(1+\beta(t))t+\tilde c+c\beta(t)})$. Then
$I_d$ is
time integral of the relative entropies of these two increment probability
measures.

Recall, for probability measures $\mu$ and $\nu$, that
the relative entropy of $\mu$ with respect to $\nu$ is defined as
%
\[
R(\mu\mmid\nu):= \cases{
\displaystyle\int\log\biggl(\frac{d\mu}{d\nu}\biggr) \,d\mu, &\quad if
$\mu\ll\nu$,\vspace*{1pt}\cr
\infty, &\quad otherwise.}
\]

%
\begin{theorem}[(Finite-dimensional LDP)]
\label{LDPmain}
The $C([0,1];\R^{d+2})$-valued sequence $\{\mathbf{X}^{n,d}\}$
satisfies an LDP with rate $n$
and convex, good rate function~$I_d$.
\end{theorem}
%
%
\begin{remark}
We now make comments on the underlying assumption (ND) and obtain the rate
function at the fixed time $t=1$.

\begin{longlist}[(A)]
\item[(A)]
The assumption (ND) specifies that the process considered is
``nondegenerate'' in some sense. (ND) does not cover some
``boundary'' cases, for instance, when $p(t)\equiv1$, the process is
deterministic in that at each time, one places a new ball in a new urn.
Also, when $\beta(t)\equiv0$, urns without a ball have no weight; and,
if in addition $p(t)\equiv0$, then all new balls are placed into urns
in the initial configuration. Although an LDP should hold in these and
other ``less degenerate'' cases, the form of the rate function may
differ in that some increments may not be possible.

On the other hand, assumption (ND) is natural with
respect to the convergence estimates needed for the proof of the lower
bound in the LDP.
However, we point out
the LDP upper bound holds without any of the boundedness assumptions on
$p(\cdot)$
and $\beta(\cdot)$ in (ND).

Formally, when $\beta(t) \equiv\infty$, this is the case of
``uniform,'' as opposed
to preferential, selection of urns. The limit
$\lim_{\beta\uparrow\infty} I_d$ corresponds to the rate function for
this type of
dynamic.

\item[(B)] One recovers the LDP at a fixed time, say $t=1$, by the
contraction principle with respect to continuous function\vspace*{1pt}
$F\dvtx C([0,1];\R^{d+2})\to\R^{d+2}$ defined by
$F(\varphi)=\varphi(1)$, so that $F(\mathbf X^{n,d})=\mathbf
X^{n,d}(1)=\frac1n \mathbf Z^{n,d}(n)$. Then Theorem~\ref{LDPmain}
implies the LDP for $\frac1n \mathbf Z^{n,d}(n)$ with rate function
given by the variational expression
$K(x) = \inf\{I_d(\varphi)\mid\varphi(0)=\mathbf{c}^d,
\varphi
(1)=x\}$ which might be evaluated numerically; cf.
\cite{Bryc-Minda-S} for calculations when $d=0$.
\end{longlist}
\end{remark}

We now extend the finite-dimen\-sional LDP results to the
infinite-dimen\-sional case ($d=\infty$). Define for
$\xi\in\prod_{i=0}^\infty C([0,1];\R)$ the function
%
%
\begin{eqnarray*}
I^{\infty}(\xi)
&=& \int_0^1 \lim_{d\rightarrow \infty} \Biggl[\bigl(1-[\dot{\xi}(t)]_0\bigr) \log
\frac{1-[\dot{\xi}(t)]_0}{p(t)+(1-p(t))\frac{\beta(t)\xi_0(t)}{(1+\beta(t))t+\tilde
c+c\beta(t)}}
\\
&&\hspace*{22pt}
{} + \sum_{i=1}^d \bigl(1-[\dot{\xi}(t)]_i\bigr) \log
\frac{1-[\dot{\xi}(t)]_i}{(1-p(t))\frac{(i+\beta(t))\xi_i(t)}{(1+\beta(t))t+\tilde
c+c\beta(t)}}\\
&&\hspace*{22pt}
{} + \Biggl(1-\sum_{i=0}^d\bigl(1- [\dot{\xi}(t)]_i\bigr)\Biggr) \log
\frac{1-\sum_{i=0}^d(1- [\dot{\xi}(t)]_i)}{(1-p(t))(1-\frac{\sum_{i=0}^d
(i+\beta(t))\xi_i(t)}{(1+\beta(t))t+\tilde c+c\beta(t)})}\Biggr]\,dt
\end{eqnarray*}
where $\xi_i(0)=c_i$, $\xi_i(t)\geq0$ is Lipschitz with constant 1,
$0\leq[{\dot\xi(t)}]_i\leq1$ for \mbox{$i\geq0$}, $\frac{d}{dt}\sum
_{i=0}^\infty
\xi_i(t) = 1$ and $\lim_d [\sum_{i=0}^d i\dot\xi_i(t) +
(d+1)(1-[\dot\xi(t)]_d)] =\sum_{i=0}^\infty(1-[\dot\xi(t)]_i)
\leq1$
for almost
all $t$, and the integral converges; otherwise $I^{\infty}(\xi
)=\infty
$. It will
turn out through a projective limit approach (cf.~\cite{Dembo-Zeitouni},
Section 4.6) that $I^\infty$ is well defined, convex and a
good rate
function, and that the integrand limit exists because the term in
square brackets is
increasing in $d$.
%
%
\begin{theorem}[(Infinite-dimensional LDP)]\label{LDPinfin}
The $\prod_{i=0}^\infty
C([0,1];\R)$-valued sequence $\{\mathbf{X}^{n,\infty}\}$ satisfies an
LDP with rate $n$
and convex, good rate function~$I^\infty$.
\end{theorem}
%
%
\begin{remark}
\label{remarkinfinity}
From the result, degree distributions, not fully supported on the
nonnegative integers, that is, when $\sum_{i\geq0} \varphi_i(t)< t + c$
or, in
other words, when the distribution specifies a positive fraction of
urns with an
infinite number of balls, cannot be achieved with finite cost in the evolution
process.
This stabilization of the ``mass'' is understood as follows.
The fraction
of urns with size larger than $A$ at time $\lfloor nt\rfloor$ is
bounded in terms of
the fraction of balls in the system: $\sum_{k\geq A} Z^n_k(\lfloor
nt\rfloor)/n \leq A^{-1}\sum_{k\geq0}
kZ^n_k(\lfloor nt\rfloor)/n \leq
A^{-1}(\lfloor nt\rfloor/n+ \tilde c^n) \leq A^{-1}(1+
2\tilde c)$ for all large $n$. Hence, for all realizations of the
process, the fraction of infinite sized urns
vanishes.

On the other hand, it seems some fraction of the total ``weight'' can
indeed be lost
in the evolution process with finite rate, that is, it may be possible
to achieve a
degree distribution at a time $0<t\leq1$ such that $\sum_{i=0}^d i\xi
_i(t) < t+\tilde c$ although
pre-limit $\sum_{i= 0}^\infty i Z^n_i(\lfloor nt\rfloor)/n = \lfloor
nt\rfloor/n +
\tilde c^n$. The interpretation is that it is possible to put a
positive fraction
of the balls into a few very large urns with finite cost, a~sort of
``condensation''
effect noticed in the limiting evolution when the selection function is
superlinear
as mentioned in Example~\ref{example1}.

The last term in the integrand of the rate function, corresponding to
the increment
${\mathbf f}^d_{d+1}$, measures the cost of choosing urns with very large
size. In the
$d\uparrow\infty$ limit, this last term may be viewed as the cost of
``escape'' of
weight from urns with bounded size, or, in other words, the cost of the
increment
``$\langle1,0,\ldots, 0,\ldots\rangle$'' which corresponds to a new
empty urn being
included and very large sized urns being incremented. Some
``condensed'' finite rate
evolutions are discussed in Example~\ref{example2}.

However, on the other hand, this type of ``weight'' loss or
``condensation'' cannot happen in the typical
evolution---see Corollary~\ref{LLN}.
\end{remark}
%
%
\begin{example}
\label{example2}
Consider the ``classical'' BA model which follows the evolution of
a random graph with preferential attachment selection function
$w(d)=d$, noted
in Example~\ref{example1} and Section~\ref{model}, which corresponds
to the
urn system when $\beta(t) \equiv1$ and $p(t)\equiv0$. Suppose that
the initial
configurations satisfy $c_i = 0$ for all $i\geq0$.

We now compute the cost of distributions in form $\xi(t) = t\gamma$ where
$\gamma=\langle\gamma_i\dvtx i\geq0\rangle$ where constants
$\gamma_i\geq0$ are such that
\[
\sum_{i\geq0}\gamma_i=1 \quad\mbox{and}\quad\sum_{i\geq
0}i\gamma_i=\sum_{i\geq0}(1-[\gamma]_i)\leq1.
\]
Since, $\xi(t)$ is linear in
$t$, calculation of the rate $I^\infty(\xi)$ simplifies considerably,
and one
evaluates the limit of the last term in the integrand of $I^\infty(\xi
)$ as the
time-independent quantity,
\begin{eqnarray*}
&&\lim_{d\uparrow\infty} \Biggl(1-\sum_{i=0}^d \bigl(1-[\dot\xi
(t)]_i\bigr)\Biggr) \log
\frac{1-\sum_{i=0}^d(1-[\dot\xi(t)]_i)}{1-(\sum
_{i=0}^d(i+1)\xi_i(t))/({2t})}\\
&&\qquad= \biggl(1-\sum_{i\geq0}i\gamma_i\biggr)\log2,
\end{eqnarray*}
which gives the cost of the ``increment $\langle1,0,\ldots, 0,\ldots
\rangle$''
when the
dynamics attaches new vertices to very large hubs or places balls into already
very large urns.

This cost is positive if $\sum_{i\geq0} i\gamma_i<1$,\vspace*{1pt} and, as
discussed in the
remark above, corresponds to the cost of forming urns/nodes with very large
size/degree in the evolution process, a ``condensation'' effect.
It follows then
%
%
\begin{equation}
\label{rateeqn}
I^\infty(\xi) = \sum_{i\geq0} (1-[\gamma]_i)\log
\frac{1-[\gamma]_i}{(i+1)\gamma_i/2}
+ \biggl(1-\sum_{i\geq
0}i\gamma
_i\biggr)\log2.
\end{equation}

In the case $\gamma_0 =1$ and $\gamma_i=0$ for $i\geq1$, one observes
$I^\infty(\xi)=\log2$, and one can associate a graph evolution to
achieve this
degree or size distribution. For instance, one may grow a ``star'' tree
configuration
where all new vertices connect to the same vertex, or all balls are put
in the same
urn. If initially, there are only two vertices with degree $1$ or two
empty urns, then the ``star'' configuration
at the $n$th step has probability $2^{-n}$ of occurring. As the degree/size
structure at time $n$ consists of $n$ leaves/empty urns and one vertex
with degree
$n$ or one urn with size $n-1$, one observes the LLN limit for the
degree/size sequence is $\xi(t) = t\gamma$, from which the rate
evaluation follows.

As discussed in Example~\ref{example1}, this ``condensed''
configuration is the limit tree with respect to superlinear
selection function $w(d) = d^r$ for $r>2$.
Moreover, as noted in the \hyperref[intro]{Introduction}, all preferential attachment
evolutions with respect to superlinear
selection function $w(d) = d^r$ for $r>1$ lead to degree distribution
$\gamma$, that is,
$E\mathcal{Z}_i(n)/n \rightarrow\gamma_i$, where $\gamma_0=1$ and
$\gamma_i=0$ for
$i\geq1$.

From formula (\ref{rateeqn}), when $\gamma$ is supported only on a
finite number of indices, one sees that
$I^\infty(\xi)<\infty$ exactly when there exists $i^*\geq0$ such that
$\gamma_{i}
>0$ for $i\leq i^*$. In particular, the ``straight road'' evolution,
leading to trees where all nodes have degree $2$, except for two leaves,
or urn configurations consisting of single ball urns except for two
empty urns, has
infinite cost: start with two vertices with degree $1$ or two empty
urns. At step $j+1$, connect a new vertex to one of the two leaves,
or add an empty urn and place a ball in one of the two empty urns in
the configuration formed at step $j$. This configuration at time $n$
has probability $1/n!$ of occurring, and in the LLN limit corresponds
to $\xi(t) = t\gamma$, where $\gamma_0=0$, $\gamma_1=1$ and $\gamma
_i =
0$ for $i\geq2$, for which $I^\infty(\xi)=\infty$.

Even when no weight escapes, that is, $\sum_{i\geq0} i\gamma_i = 1$, it
may be noted
that deviations to nonpower law urn size paths $\xi$ are possible
with finite
rate. For instance, when $\gamma_i = 2^{-(i+1)}$ for $i\geq0$,
$I^\infty(\xi) =
-\sum_{i\geq0} \frac{1}{2^{i+1}}\log\frac{i+1}{2}$. When $\gamma
_i =
q(i+1)$ for
$i\geq0$ and $q$ in form of the stretched exponential in
(\ref{stretched}), the LLN limit for the degree distribution with
respect to sublinear selection preferential attachment, a calculation
verifies that $\sum_{i\geq0}i\gamma_i=1$ and also $I^\infty(\xi
)<\infty$.
\end{example}


We now turn to the LLN behavior
which corresponds to the ``zero-cost'' trajectory.
Consider the system of ODEs for $\varphi^d = \varphi$, with initial
condition $\varphi(0)=\mathbf c^d$:
%
%
\begin{eqnarray}
\label{ODEforLLN}
\dot\varphi_0(t)
& =& 1-p(t)-\bigl(1-p(t)\bigr)\frac{\beta(t)\varphi_0(t)}{(1+\beta
(t))t+\tilde
c+c\beta(t)},\nonumber\\
\dot\varphi_1(t)
& =& p(t)+\bigl(1-p(t)\bigr)\frac{\beta(t)\varphi_0(t)}{(1+\beta(t))t+\tilde
c+c\beta(t)}\nonumber\\
&&{} - \bigl(1-p(t)\bigr)\frac{(1+\beta(t))\varphi_1(t)}{(1+\beta
(t))t+\tilde c+c\beta(t)},\nonumber\\[-8pt]\\[-8pt]
\dot\varphi_i(t)
& =& \bigl(1-p(t)\bigr)\frac{(i-1+\beta(t))\varphi_{i-1}(t)}{(1+\beta
(t))t+\tilde c+c\beta(t)}\nonumber\\
&&{} - \bigl(1-p(t)\bigr)\frac{(i+\beta(t))\varphi_i(t)}{(1+\beta
(t))t+\tilde c+c\beta(t)}\qquad \mbox{for }2\leq i\leq d,\nonumber\\
\dot\varphi_{d+1}(t) & =& 1 - \sum_{i=0}^d \dot\varphi
_i(t).\nonumber
\end{eqnarray}
Recall that a ``Carath\'eodory'' solution is an absolutely continuous
function satisfying the ODEs a.a. $t$, and the initial condition, or
equivalently a function satisfying the integral equation associated to
the ODEs.
One can readily integrate ODEs (\ref{ODEforLLN}), and find a
Carath\'{e}odory solution $\zeta^d(t) = \langle\zeta_0(t),\ldots, \zeta
_d(t),\bar
\zeta
_{d+1}(t)\rangle$ [see formula (\ref{zero-cost})], which is unique
from the
following theorem. One extends to ``$d=\infty$'' setting by defining
\[
\zeta^\infty(t):=
\langle\zeta_0(t),\ldots,\zeta_d(t),\ldots\rangle.
\]

We now state a LLN for $\mathbf X^{n,d}$ and $\mathbf X^{n,\infty}$ as
a consequence of the LDP upper bound. As remarked in the
\hyperref[intro]{Introduction}, this LLN may also be obtained by rate equation
formulations as in~\cite{KR} and~\cite{Du06}, Chapter 4.
%
%
\begin{corollary}[(LLN)]\label{LLN}
For $d\geq0$, $\zeta^d$ is the unique Carath\'eodory solution to ODEs
(\ref{ODEforLLN}) with the initial condition $\varphi(0)=\mathbf
c^d$, and also $I_d(\zeta^d)=0$. Then, in the sup topology on
$C([0,1];\R^{d+2})$,
$\mathbf{X}^{n,d}(\cdot) \to\zeta^d(\cdot)$ a.s.

As a consequence, we have in the product topology on
$\prod_{i=0}^\infty C([0,1]; \R)$ that $\mathbf{X}^{n,\infty}(\cdot)
\rightarrow\zeta^\infty(\cdot)$.
Moreover,
$\sum_{i=0}^\infty\zeta_i(t) = t + c$ and $\sum_{i=0}^\infty i\zeta
_i(t) = t + \tilde c$, and hence no ``weight'' is lost in the LLN limit.
\end{corollary}
%
%
\begin{remark}
\label{LLNrmk}
The last equality, $\sum_{i\geq0} i\zeta_i(t) = t+\tilde c$, requires
the condition in (LIM) that the
initial scaled degree distribution is not ``condensed,'' that is,
$\tilde c =
\lim_{n\uparrow\infty} \tilde c^n = \sum_{i\geq0} i c_i$. When
the initial distribution is ``condensed,'' that is, a strict Fatou
limit $\tilde c =
\lim_{n\uparrow\infty} \tilde c^n > \sum_{i\geq0} ic_i$ occurs,
the large
deviation results Theorems~\ref{LDPmain},~\ref{LDPinfin} and Corollary
\ref{LLN} (except for the last equality) still hold with the same
notation and proofs. However, one can show by similar arguments as
for the proof of the last equality in Corollary~\ref{LLN} that
the LLN trajectory $\zeta^\infty$ will now be
``condensed,'' that is, $s(t)=\sum_{i\geq0} i\zeta_i(t) < t+
\tilde c$ for $t\geq0$. Moreover, for a constant $C=C(c,\tilde c,
p_0, \beta_1,\beta_0)>0$, one can see for all large $t$ that
\[
C\biggl(\tilde c - \sum_{i\geq0}ic_i\biggr)t^{(1-p_0)/(1+\beta_1)}
\leq
t+\tilde c - s(t) \leq C^{-1}\biggl(\tilde c - \sum_{i\geq
0}ic_i\biggr)
t^{1/(1+\beta_0)}.
\]
\end{remark}

%
\begin{figure}

\includegraphics{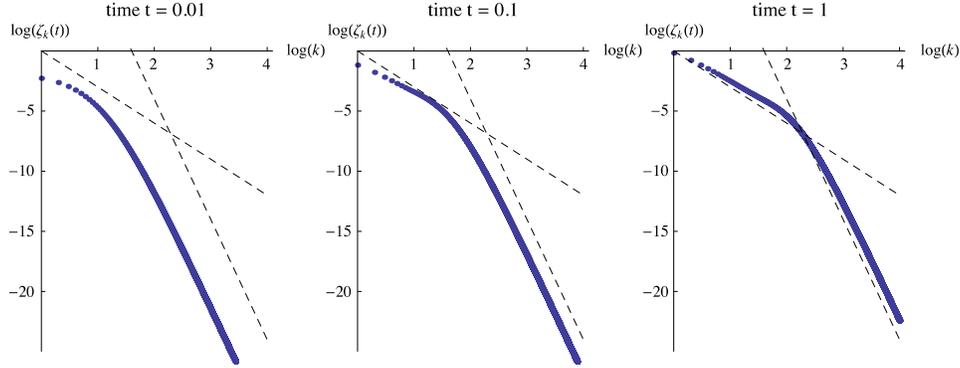}

\caption{The thick curves are the (numerical) LLN ODE paths at times
$t=0.01,0.1,1$
with $p(t)\equiv0$, $\beta(t)=8$ for $t<0.01$, $1$ for $t\geq0.01$
and $c_k\equiv0$. Dashed lines are straight lines with slope $-3$ and
$-10$. The plots use log--log scales.}\label{fig1}
\end{figure}

We now consider the ``scale-freeness'' of $\zeta^\infty$. Although it
seems difficult to control each $\zeta_i$, nevertheless $\zeta^\infty$
has ``power law'' behavior, in terms of bounds on $[\zeta^\infty]_i$.
In general, it appears $\zeta^\infty$ can interpolate between the
bounds (cf. Figure~\ref{fig1}; as a curiosity, we note a figure with a similar
``bend'' is found in~\cite{Gjoka-Kurant-Butts-Markopoulou} with respect
to Facebook social network data).
%
%
\begin{corollary}[(Power law)]\label{powerlaw}
Assume
$0\leq p_{\min}\leq p(\cdot)\leq p_0=:p_{\max}<1$, and $0< \beta
_0=:\beta_{\min}\leq\beta(\cdot)\leq\beta_{\max}:=\beta_1$. Then,
$\zeta^\infty$ is bounded between two power laws:

For small configurations, for example, $c_k \equiv0$, we
have, for $i\geq0$ and $t\geq0$,
\[
[\eta']_i t \leq[\zeta^\infty(t)]_i \leq[\eta]_i t.
\]
For large configurations, for example, $c_k>0$ for some $k\geq0$, we
have, for $i\geq0$,
\[
[\eta']_i \bigl(t+o(1)\bigr) \leq[\zeta^\infty(t)]_i \leq
[\eta]_i \bigl(t+o(1)\bigr) \qquad\mbox{as } t\uparrow\infty.
\]
Here, with respect to positive constants $C, C'$ depending on $p$ and
$\beta$,
\[
\eta'_i:= \frac{C'}{i^{1+({1+\beta_{\min}})/({1-p_{\min
}})}}\bigl(1+o(1)\bigr)
\]
and
\[
\eta_i:= \frac{C}{i^{1+({1+\beta_{\max}})/({1-p_{\max}})}}\bigl(1+o(1)\bigr).
\]
\end{corollary}

The outline of the paper is that in Sections~\ref{sectionLDP} and
\ref{sectionLDPinfin}, we prove the finite and
infinite-dimensional LDPs, Theorems~\ref{LDPmain} and~\ref{LDPinfin}.
In Section~\ref{sectionLLN}, we prove the law of large numbers
(Corollary~\ref{LLN}). Finally, in Section~\ref{sectionpowerlaw}, we
discuss power-law behavior (Corollary~\ref{powerlaw}).

\section{\texorpdfstring{Proof of Theorem \protect\ref{LDPmain}}{Proof of Theorem 1.2}}\label{sectionLDP}

We follow the method and notation of Dupuis and Ellis in
\cite{Dupuis-Ellis}; see also~\cite{Zhang-Dupuis-08}. Some steps are similar
to those in
\cite{Bryc-Minda-S} where the ``leaves'' in a more simplified graph
scheme are considered. However, as many things differ in our model, in
the upper bound, and especially the lower bound proof, we present the
full argument.

We now fix $0\leq d<\infty$ and equip $\R^{d+2}$ with the $L_1$-norm
denoted by $|\cdot|$.
Recall, from assumption (LIM),
\[
\mathbf c^{n,d}=(c^{n}_0,\ldots,c^n_d,\bar c^{n,d}):=\frac1n \mathbf
Z^{n,d}(0)
\to\mathbf c^d,
\]
where $\bar c^{n,d} = \sum_{i\geq d+1}c_i^n$.
Denote
\[
\vec\xi(n,t):= (p_n(t),\beta_n(t),\sigma_n(t)),
\]
where
\begin{eqnarray*}
p_n(t)&:=&p({\lfloor nt\rfloor}/n),\qquad \beta_n(t):= \beta
({\lfloor
nt\rfloor}/n),\\
\sigma_n(t)&:=&\frac1n s^n(\lfloor nt \rfloor)=\bigl(1+\beta_n(t)\bigr)\frac
{\lfloor nt\rfloor}n + \tilde c^{n} + c^{n}\beta_n(t).
\end{eqnarray*}
Let
\begin{eqnarray*}
\sigma(t)&:=& \bigl(1+\beta(t)\bigr)t+\tilde c + c\beta(t),\\
\vec\xi(t)&:=& (p(t),\beta(t),\sigma(t)).
\end{eqnarray*}
We note that, as $n\to\infty$, and $p(t)$ and $\beta(t)$ are piecewise
continuous,
\[
\vec\xi(n,t)\to\vec\xi(t) \qquad\mbox{for almost all }t\in[0,1].
\]

In the remainder of the section, when the context is clear, we often
drop the superscript $d$ to save on notation.
Recall
\[
\mathbf{X}^n(j):=\frac1n \mathbf Z^{n,d}(j),
\]
$\mathbf X^n(0)=\mathbf c^{n,d}$ and $\mathbf{X}^n(j+1)=\mathbf
{X}^n(j)+\frac1n y_{\mathbf{X}^n(j)}^n(j)$,
where
\[
y_{\mathbf{x}}^n(j) \mbox{ has distribution }\rho_{\vec\xi
(n,j/n), \mathbf{x}}.
\]
Here, for
$\mathbf{x}=\langle x_0,\ldots, x_d, x_{d+1}\rangle\in\R^{d+2}$ such that
$x_i\geq
0$ for $0\leq i\leq d+1$, numbers $p'\in[0,1]$ and $\beta',\sigma
'\geq
0$ such that $\sum_{i=0}^{d+1} (i+\beta')x_i \leq\sigma'$, and
\mbox{$A\subset\R^{d+2}$},
\begin{eqnarray*}
\rho_{(p',\beta',\sigma'), \mathbf{x}}(A)&:=&
\biggl(p' +(1-p') \frac{\beta' x_0}{\sigma'}\biggr)\delta_{\mathbf
f_0}(A)\\
&&{} + \sum_{i=1}^d(1-p')\frac{(i+\beta')x_i}{\sigma'}\delta
_{\mathbf f_i}(A)\\
&&{} + (1-p')\biggl(1-\frac{\sum_{i=0}^d (i+\beta
')x_i}{\sigma
'}\biggr)\delta_{\mathbf f_{d+1}}(A).
\end{eqnarray*}
We note when $\sigma'=0$ and $\mathbf{x}=\langle0,\ldots,0\rangle$,
by convention
$0/0=0$ and
\[
\rho_{(p',\beta',0),\mathbf{x}}(A):=p'\delta_{\mathbf
{f}_0}(A) + (1-p')\delta_{\mathbf{f}_{d+1}}(A).
\]

From (\ref{increments}) and (LIM), for $A>0$, the paths $\mathbf{X}^n(t)=
\mathbf{X}^{n,d}(t)$, for all large~$n$,
belong to
%
%
\begin{eqnarray}\label{spaceforpaths}
\Gamma_{d,A} &:=& \Biggl\{\varphi\in C([0,1];\R^{d+2}) |
|\varphi
(0)-\mathbf c^d| \leq A, \varphi_i \mbox{ is Lipschitz}\nonumber
\\
&&\hspace*{8pt}\mbox{with bound }1, 0\leq[\dot\varphi(t)]_i\leq1
\mbox{ for } 0\leq
i\leq d+1\mbox{, and}\\
&&\hspace*{8pt}\sum_{i=0}^{d+1}\dot\varphi_i(t)= 1, \sum_{i=0}^{d+1} i\dot
\varphi
_i(t)= \sum_{i=0}^d \bigl(1-[\dot\varphi(t)]_i\bigr)\leq1 \mbox{ for a.a.
}t\Biggr\}.\nonumber
\end{eqnarray}
Here, we equip $C([0,1];\R^{d+2})$ with the supremum norm.

Let $h\dvtx C([0,1];\R^{d+2}) \to\R$ be a bounded
continuous function. Let also
\[
W^n:=-\frac1 n \log E\{\exp[-n h(\mathbf{X}^n)]\}.
\]
To prove Theorem~\ref{LDPmain}, we need to establish Laplace principle
upper and lower
bounds (cf.~\cite{Dupuis-Ellis}, Section 1.2), namely upper bound
\[
\liminf_{n\to\infty} W^n \geq\inf_{\varphi\in C([0,1];\R^{d+2})}
\{I_d(\varphi) + h(\varphi)\}
\]
for a good rate function $I_d$, and lower bound
\[
\limsup_{n\to\infty} W^n \leq\inf_{\varphi\in C([0,1];\R^{d+2})}
\{I_d(\varphi) + h(\varphi)\}.
\]

Given $\mathbf X^n(0) = \mathbf c^{n,d}$, define, for $1\leq j\leq n$, that
\begin{eqnarray*}
&& W^n(j,\{\mathbf x_{1},\ldots, \mathbf x_j \})\\
&&\qquad:= -\frac1 n \log E\{\exp[-n h(\mathbf X^n )]
\mid\mathbf X^n(1) =\mathbf x_{1},\ldots, \mathbf{X}^n(j)=\mathbf
x_j \}
\end{eqnarray*}
and
\[
W^n:= W^n(0,\varnothing) = -\frac1 n \log E\{\exp[-n h(\mathbf X^n
)]\}.
\]
The Dupuis--Ellis method stems from the following discussion. From the
Markov property,
for $1\leq j\leq n-1$,
\begin{eqnarray*}
&&e^{-nW^n(j,\{\mathbf x_{1},\ldots, \mathbf x_j \})} \\
&&\qquad= E\bigl\{e^{-n
h(\mathbf X^n )}\mid\mathbf X^n(1) =\mathbf x_{1}, \ldots, \mathbf
{X}^n(j)=\mathbf x_j \bigr\}\\
&&\qquad = E\bigl\{ E\bigl\{e^{-n h(\mathbf X^n )}\mid\mathbf X^n(1),\ldots,
\mathbf X^n(j+1)\bigr\}\mid\mathbf X^n(1)=\mathbf x_{1},\ldots, \mathbf
{X}^n(j)=\mathbf x_j\bigr\}\\
&&\qquad = E\bigl\{ e^{-nW^n(j+1,\{\mathbf X^n(1),\ldots, \mathbf{X}^n(j),
\mathbf X^n(j+1)\})}\mid\mathbf
X^n(1)=\mathbf x_{1},\ldots, \mathbf{X}^n(j)=\mathbf x_j\bigr\}\\
&&\qquad = \int_{\R^{d+2}} e^{-nW^n(j+1,\{\mathbf x_{1},\ldots,
\mathbf
x_j, \mathbf x_j + \mathbf{y}/n\})} \rho_{\vec\xi(n,j/n),
\mathbf x_j}(d\mathbf{y}).
\end{eqnarray*}
Recall the definition of relative entropy near Theorem~\ref{LDPmain}.
Then, by the variational formula for relative entropy (cf.
\cite{Dupuis-Ellis}, Proposition 1.4.2), for $1\leq j\leq n-1$,
\begin{eqnarray*}
&&W^n(j,\{\mathbf x_{1},\ldots, \mathbf x_j \})\\
&&\qquad = -\frac1 n \log\int_{\R^{d+2}} e^{-nW^n(j+1,\{\mathbf
x_{1},\ldots, \mathbf x_j, \mathbf x_j + \mathbf{y}/n\})}
\rho
_{\vec\xi(n,j/n), \mathbf x_j}(d\mathbf{y})\\
&&\qquad = \inf_{\mu}\biggl\{\frac1 n R\bigl(\mu\mmid\rho_{\vec\xi(n,j/n),
\mathbf x_j}\bigr)\\
&&\hspace*{16.3pt}\qquad\quad{} + \int_{\R^{d+2}} W^n\biggl(j+1,\biggl\{\mathbf x_{1},\ldots,
\mathbf x_j, \mathbf x_j
+ \frac1 n \mathbf{y}\biggr\}\biggr) \mu(d\mathbf{y})\biggr\}.
\end{eqnarray*}
We also have a terminal condition $W^n(n,\{\mathbf x_{1},\ldots
,\mathbf
x_n\}) = h(\mathbf x(\cdot))$, where $\mathbf x(\cdot)$
is the linear interpolated path connecting $\{(j/n,\mathbf x_j)\}
_{0\leq j\leq n}$.

We may understand these dynamic programming equations and terminal
conditions in terms of a particular stochastic control problem.
Define:
\begin{longlist}
\item$\mathcal L_j=(\R^{d+2})^{j}$, the state space on which
$W^n(j,\cdot)$ is defined;
\item$\mathcal U= \mathcal P(\R^{d+2})$, where $\mathcal P(B)$ is the
space of probabilities on $B$, is the control space on which the
infimum is taken;
\item for $j={0},\ldots,n-1$, ``control'' $v_j^n(d\mathbf
{y})=v_j^n(d\mathbf{y}|\mathbf{x}_{0},\ldots,\mathbf{x}_j)$
which is a stochastic kernel on $\R^{d+2}$ given $(\R^{d+2})^{j}$;
\item$\{\bar{\mathbf X}^n(j); 0\leq j\leq n\}$, the ``controlled''
process which is the adapted path satisfying
$\bar{\mathbf X}^n(0)=\mathbf c^{n,d}$ and $\bar{\mathbf X}^n(j+1) =
\bar{\mathbf X}^n(j) + \frac1 n \bar{\mathbf Y}^n(j)$ for $0\leq
j\leq n-1$,
where $\bar{\mathbf Y}^n(j)$, given $(\bar{\mathbf X}^n(0), \ldots,
\bar{\mathbf X}^n(j))$,
has distribution $v_j^n(\cdot)$ [i.e., $\bar P\{\bar{\mathbf
Y}^n(j)\in d\mathbf{y}\mid\bar{\mathbf X}^n(0), \ldots, \bar
{\mathbf
X}^n(j)\}
:= v_j^n(d\mathbf{y}\mid\bar{\mathbf X}^n(0), \ldots, \bar{\mathbf
X}^n(j))$] and
$\bar{\mathbf X}^n(\cdot)$ is the piecewise linear interpolation of
$(j/n,\bar{\mathbf X}^n(j))$;\vspace*{1pt}%
\item``running costs'' $C_j(v)=\frac1 n R(v\mmid\rho)$ for $v\in
\mathcal P(\R^{d+2})$; and\vspace*{1pt}
\item``terminal cost'' equals to the function $h$.%
\end{longlist}

Also, define, for $0\leq j\leq n-1$, the minimal cost function
%
\begin{eqnarray*}
&&V^n(j,\{\mathbf{x}_1,\ldots, \mathbf{x}_j\})\\[-2pt]
&&\qquad =
\inf_{\{v_i^n\}} \bar E_{j, \mathbf{x}_1,\ldots, \mathbf{x}_j}
\Biggl\{ \frac1 n \sum_{i=j}^{n-1} R\bigl(v_i^n(\cdot) \mmid\rho
_{\vec\xi
(n,i/n), \bar{\mathbf X}^n(i)}\bigr) + h(\bar{\mathbf X}^n(\cdot
))\Biggr\},
\end{eqnarray*}
where $v_i^n(\cdot)=v_i^n(\cdot\mid\bar{\mathbf X}^n(0), \ldots,
\bar
{\mathbf X}^n(i))$, and
the infimum is taken over all control
sequences $\{v_i^n\}$. Here,
$\bar E_{j, \mathbf{x}_1,\ldots, \mathbf{x}_j}$ denotes
expectation, with respect to the adapted process $\bar{\mathbf
{X}}^n(\cdot)$ associated to $\{v_i^n\}$, conditioned
on $\bar{\mathbf X}^n(1)=\mathbf{x}_1,\ldots, \bar{\mathbf X}^n(j) =
\mathbf{x}_j$.
The boundary conditions are
$V^n(n,\{\mathbf{x}_1,\ldots, \mathbf{x}_n\}) =
h(\mathbf{x}(\cdot))$ and
%
%
\begin{equation}\label{eqrepformulafornonleaves}\qquad
V^n:= V^n(0,\varnothing)= \inf_{\{v_j^n\}} \bar E \Biggl\{ \frac1
n \sum_{j=0}^{n-1}
R\bigl(v_j^n(\cdot) \mmid\rho_{\vec\xi(n,j/n), {\bar{\mathbf
X}^n(j)}}\bigr)
+ h({\bar{\mathbf X}^n(\cdot)})\Biggr\}.
\end{equation}

It turns out that $\{V^n(j,\{\mathbf{x}_1,\ldots, \mathbf{x}_j\})\dvtx
0\leq j\leq n\}$ also satisfies the dynamic programming equations and
terminal condition, and since these equations have unique solutions
(cf.~\cite{Dupuis-Ellis}, Section 3.2), we may conclude by
\cite{Dupuis-Ellis}, Corollary 5.2.1, that
%
\[
W^n = -\frac1 n \log E\{\exp[-n h(\bar{\mathbf{X}}^n(\cdot))]\} = V^n.
\]


\subsection{Upper bound}

To prove the upper bound, it will be helpful to put the
controls $\{v_j^n\}$ into continuous-time paths. Let
$v^n(d\mathbf{y}|t):=v_j^n(d\mathbf{y})$ for $t\in[j/n, (j+1)/n)$,
$j=0,\ldots,n-1$, and $v^n(d\mathbf{y}|1):=v_{n-1}^n$. Define
\[
v^n(A\times B):=\int_B v^n(A|t)\,dt
\]
for Borel $A\subset\R^{d+2}$ and $B\subset[0,1]$. Also define the
piecewise constant path ${\tilde{\mathbf X}^n}(t):={\bar{\mathbf
X}^n(j)}$ for $t\in[j/n, (j+1)/n)$, $0\leq j\leq n-1$, and
${\tilde{\mathbf X}^n}(1):={\bar{\mathbf X}^n(n-1)}$. Then
\[
W^n = V^n = \inf_{\{v_j^n\}} \bar E \biggl\{ \int_0^1
R\bigl(v^n(\cdot\mid t) \mmid\rho_{\vec\xi(n,t),{\tilde{\mathbf
X}^n}(t)}\bigr) \,dt
+ h({\bar{\mathbf X}^n})\biggr\}.
\]
Given $\rho_{\vec\xi,\mathbf{x}}$ is supported on $K:=\{\mathbf
f_0,\mathbf f_1,\ldots,\mathbf{f}_{d+1}\}$,
if ${\{v_j^n\}}$ is not supported on $K$, then
$R(v^n\mmid\rho_{\vec\xi,\mathbf{x}})=\infty$. Since $|V^n|\leq
\|h\|_\infty<\infty$ and $K\subset\R^{d+2}$ is compact, for each $n$,
there is $\{v_j^n\}$ supported
on $K$ and corresponding $v^n(d\mathbf{y}\times
dt)=v^n(d\mathbf{y}\mid t)\times dt$ such that, for $\varepsilon>0$,
%
%
\begin{equation}\label{eqUB}
W^n + \varepsilon= V^n + \varepsilon\geq\bar E \biggl\{ \int_0^1
R\bigl(v^n(\cdot\mid t) \mmid\rho_{\vec\xi(n,t), \tilde{\mathbf
X}^n(t)}\bigr) \,dt
+ h(\bar{\mathbf X}^n)\biggr\}.%
\end{equation}

Recall that $\bar{\mathbf X}^{n}(\cdot)$ takes values in
$\Gamma_{d,A}$. Since $\Gamma_{d,A}$ is compact, by applications of
the Ascoli--Arzel\'a theorem, and $\{v_j^n\}$ is tight, by
Prokhorov's theorem, given any subsequence
of $\{v^n, \bar{\mathbf X}^n\}$, there is a further subsubsequence, a
probability %
space $(\bar\Omega, \bar{\mathcal F}, \bar P)$,
a stochastic kernel $v$ on $K\times[0,1]$ given $\bar\Omega$ and a
random variable $\bar{\mathbf X}$
mapping $\bar\Omega$ into\vspace*{1pt} $\Gamma_{d, A}$ such that the
subsubsequence converges in distribution to $(v,\bar{\mathbf X})$.
In particular, since
$\bar{\mathbf X}^{n}(0) = \mathbf c^{n,d}\to\mathbf c^d$ as $n\to
\infty$, we have $\bar{\mathbf X}$ [cf. (\ref{spaceforpaths})]
belongs to
%
\[
\Gamma_d:= \Gamma_{d,0}\qquad \mbox{those functions such that
} \varphi(0)= \mathbf c^d.
\]
Then,~\cite{Dupuis-Ellis}, Lemma 3.3.1, shows that $v$
is a subsequential weak limit of $v^n$, and there exists a
stochastic kernel $v(dy\mid t,\omega)$ on $K$ given $[0,1]\times
\bar\Omega$ such that $\bar P$-a.s. for $\omega\in\bar\Omega$,
\[
v(A\times B\mid\omega) = \int_B v(A\mid t,\omega) \,dt.
\]
Now, the same
proof given for~\cite{Dupuis-Ellis}, Lemma 5.3.5, shows that $(v^n,
\bar{\mathbf X}^n, \tilde{\mathbf X}^n)$ has a subsequential weak
limit $(v,\bar{\mathbf X}, \bar{\mathbf X})$, where the last coordinate
is with respect to Skorokhod space $D([0,1];\R^{d+2})$, and $\bar P$-a.s.
for $t\in[0,1]$
%
\begin{eqnarray*}
\bar{\mathbf X}(t) &=& \int_{\R^{d+2}\times[0,t]} \mathbf{y}
v(d\mathbf{y}\times
ds) = \int_0^t \biggl(\int_K \mathbf{y} v(d\mathbf{y}\mid s)
\biggr)\,ds, \\[-2pt]
\dot{\bar{\mathbf X}}(t) &=& \int_K
\mathbf{y} v(d\mathbf{y}\mid t).
\end{eqnarray*}
By Skorokhod's representation
theorem, we may take that
$(v^n,\bar{\mathbf X}^{n}, \tilde{\mathbf X}^{n})$ converges to
$(v,\bar{\mathbf X}, \bar{\mathbf X})$ a.s. In particular,
$\bar{\mathbf X}^{n} \to\bar{\mathbf X}$ uniformly a.s., and as
$\bar{\mathbf X}$ is continuous, it follows that also $\tilde{\mathbf
X}^n \to\bar{\mathbf X}$ uniformly a.s.; cf.
\cite{Dupuis-Ellis}, Theorem~A.6.5.

Let $\lambda$ denote Lebesgue measure on $[0,1]$ and $\rho\times
\lambda
$ product measure on $K\times[0,1]$. Then
\cite{Dupuis-Ellis}, Lemma 1.4.3(f), yields
\[
\int_0^1 R\bigl(v^n(\cdot\mid t)\mmid\rho_{\vec\xi(n,t), \tilde
{\mathbf
X}^n(t)} \bigr) \,dt
= R\bigl( v^n(\cdot\mid t)\times\lambda(dt)\mmid\rho_{\vec\xi(n,t),
\tilde{\mathbf X}^n(t)} \times\lambda(dt)\bigr).
\]

We now evaluate the limit inferior of $W^n$ using formula
(\ref{eqUB}), along a subsequence as above:
\begin{eqnarray*}
\liminf_{n\to\infty} V^n +\eps
&\geq& \liminf_{n\to\infty} \bar E\biggl\{
\int_{0}^1 R\bigl(v^n(\cdot|t )\mmid\rho_{\vec\xi(n,t), \tilde
{\mathbf
X}^n(t)} \bigr) \,dt
+ h(\bar{\mathbf X}^n)
\biggr\}\\
&=& \liminf_{n\to\infty} \bar E\bigl\{
R\bigl( v^n(\cdot\mid t)\times\lambda(dt)\mmid\rho_{\vec\xi(n,t),
\tilde
{\mathbf X}^n(t)} \times\lambda(dt)\bigr)
+ h(\bar{\mathbf X}^n)
\bigr\}\\
&\geq&\bar E\bigl\{
R\bigl( v(\cdot\mid t)\times\lambda(dt)\mmid\rho_{\vec\xi(t), \bar
{\mathbf X}(t)} \times\lambda(dt)\bigr)
+ h(\bar{\mathbf X})
\bigr\}\\
&=& \bar E\biggl\{\int_{0}^1 R\bigl(v(\cdot|t )\mmid\rho
_{\vec
\xi(t), \bar{\mathbf X}(t)} \bigr) \,dt
+ h(\bar{\mathbf X})\biggr\}.%
\end{eqnarray*}
Note that we used Fatou's lemma in the second inequality, observing \mbox{(i)--(iv)}.

\begin{longlist}
\item$v^n(d\mathbf{y}|dt)\times\lambda(dt) \to v(d\mathbf
{y}|dt)\times\lambda(dt)$ a.s. as $v^n \Rightarrow v$ a.s.;\vspace*{1pt}
\item$\rho_{\vec\xi(n,t), \tilde{\mathbf X}^{n}(t)} \Rightarrow
\rho
_{\vec\xi(t), \bar{\mathbf X}(t)}$
as $\vec\xi(n,t) \to\vec\xi(t)$ a.a.\vspace*{1pt} $t\in[0,1]$, and $\tilde
{\mathbf X}^{n}(t) \to\bar{\mathbf X}(t)$ uniformly on $[0,1]$
a.s.; %
\item$\liminf_{n\to\infty} R(v^n(d\mathbf{y}|dt)\times
\lambda
(dt) \mmid\rho_{\vec\xi(n,t), \tilde{\mathbf X}^{n}(t)}
\times\lambda(dt)) \geq R(v(d\mathbf{y}|dt)\times\lambda
(dt) \mmid\rho_{\vec\xi(t),\bar{\mathbf X}(t)} \times\lambda
(dt)
)$ a.s. as
$R$ is lower semi-continuous;\vspace*{1pt}
\item$h(\bar{\mathbf X}^{n}) \to h(\bar{\mathbf X})$ a.s. as $h$ is
continuous and $\bar{\mathbf X}^{n} \to\bar{\mathbf X}$ uniformly on
$[0,1]$ a.s.
\end{longlist}

By~\cite{Dupuis-Ellis}, Lemma 3.3.3(c),
\[
R\bigl(v(\cdot|t)\mmid\rho_{\vec\xi(t), \bar{\mathbf X}(t)}\bigr)
\geq L\biggl(\vec\xi(t), \bar{\mathbf X}(t), \int_K \mathbf{z}
v(d\mathbf{z}|t)\biggr),
\]
where
\begin{eqnarray*}
L(\vec\xi(t),\mathbf{x},\mathbf{y})
:\!&=& \sup\biggl\{\langle\bolds{\theta},
\mathbf{y}\rangle
- \log\int_K \exp\langle\bolds{\theta}, \mathbf{z}\rangle\rho
_{\vec
\xi(t),\mathbf{x}}(d\mathbf{z}) \Big| \bolds{\theta}\in
\R
^{d+2}\biggr\}\\
&=& \inf\biggl\{R\bigl(\nu(\cdot|t) \mmid\rho_{\vec\xi(t),\mathbf
{x}}\bigr)
| \nu(\cdot|t)\in
\mathcal P(K), \int_K \mathbf{z}\nu(d\mathbf{z}|t)=\mathbf y
\biggr\}.
\end{eqnarray*}
We note, in this definition,
the infimum is attained at
some $\nu_0\in\mathcal{P}(K)$ as the relative entropy
is convex and lower semicontinuous; cf.~\cite{Dupuis-Ellis}, Lem\-ma~1.4.3(b).
Since $\int\mathbf z v(d\mathbf z|t)=\dot{\bar{\mathbf X}}(t)$,
we have
\[
\liminf_{n\to\infty} V^n \geq\bar E\biggl\{\int_{0}^1
L( \vec\xi(t), \bar{\mathbf X}(t),\dot{\bar{\mathbf X}}(t)) \,dt +
h(\bar{\mathbf X}) \biggr\}.
\]
As $\bar{\mathbf X} \in\Gamma_d$, we have
\[
\liminf_{n\to\infty} V^n \geq\inf_{\varphi\in\Gamma_d}
\biggl\{
\int_{0}^1
L(\vec\xi(t),\varphi(t),\dot{\varphi}(t)) \,dt +
h(\varphi)\biggr\}.
\]

For $\varphi\in\Gamma_d$, we can evaluate a unique minimizer $\nu
_0(\cdot|t)$ in the definition of
$L(\vec\xi(t),\varphi(t), \dot{\varphi}(t))$: recall that
$[\dot{\varphi}(t)]_i:=\sum_{l=0}^i \dot{\varphi}_l(t)$. Then, as
$\sum
_{i=0}^{d+1} \mathbf{f}_i \nu_0(\mathbf{f}_i|t) = \langle\dot\varphi
_0(t),\ldots,\dot\varphi_{d+1}(t)\rangle$, a calculation gives
%
%
\begin{equation}\label{nuforLbounded}
\nu_0(\dot{\varphi}(t)|t) = \sum_{i=0}^d \bigl(1-[\dot{\varphi
}(t)]_i\bigr)\delta_{\mathbf f_i} +
\Biggl({\sum_{i=0}^d [\dot{\varphi}(t)]_i-d}\Biggr)
\delta_{\mathbf f_{d+1}}.
\end{equation}
Hence,
%
%
%
\begin{eqnarray}\label{Lboundedcase}
&&L(\vec \xi(t),\varphi(t), \dot{\varphi}(t))\nonumber\\
&&\qquad= R\bigl(\nu_0(\dot\varphi(t)|t)\Vert \rho_{\vec \xi(t),\varphi(t)}\bigr)\nonumber\\
&&\qquad= \bigl(1-[\dot{\varphi}(t)]_0\bigr) \log \frac{1-[\dot{\varphi}(t)]_0}{{p(t)+(1-p(t))
\frac{\beta(t)\varphi_0(t)}{(1+\beta(t))t+\tilde c+c\beta(t)}}}\\
&&\qquad\quad{} + \sum_{i=1}^d \bigl(1-[\dot{\varphi}(t)]_i\bigr) \log
\frac{1-[\dot{\varphi}(t)]_i}{(1-p(t))\frac{(i+\beta(t))\varphi_i(t)}{(1+\beta(t))t+\tilde c+c\beta(t)}}\nonumber\\
&&\qquad\quad{} + \Biggl(1-\sum_{i=0}^d\bigl(1- [\dot{\varphi}(t)]_i\bigr)\Biggr) \log
\frac{1-\sum_{i=0}^d(1- [\dot{\varphi}(t)]_i)}{(1-p(t))
(1-\frac{\sum_{i=0}^d (i+\beta(t))\varphi_i(t)}{(1+\beta(t))t+\tilde
c+c\beta(t)})},\nonumber
\end{eqnarray}
interpreted under our conventions (\ref{convention}).

Finally, define
\[
I_d(\varphi):= \int_0^1 L( \vec\xi(t),\varphi(t), \dot{\varphi
}(t)) \,dt,
\]
when $\varphi\in\Gamma_d$, and $I_d(\varphi)=\infty$ otherwise.
Since $L$ is convex, $I_d$ is convex. Also $I_d$ has
compact level sets by the proof of
\cite{Dupuis-Ellis}, Proposition
6.2.4, and so is a good rate function.
Hence, the Laplace principle upper bound holds with respect to $I_d$.

We will need the following result for the proof of the lower bound in
the next section.
%
%
\begin{lemma}
\label{Iforlinear}
Let $\ell(t)=\mathbf e t + \mathbf c^d$ be a linear function, where
$\mathbf e=(e_0,e_1,\ldots,\break e_{d+1})$ is such that $e_i>0$
for $i\geq0$,
$\sum_{i=0}^{d+1} e_i =1$, and $\sum_{i=0}^{d+1} ie_i\leq1$.
Then, $I_d(\ell(t))<\infty$.
\end{lemma}
\begin{pf}
Noting $\sum_{i=0}^d (1-[\mathbf e]_i) = \sum_{i=0}^{d+1}ie_i\leq1$,
explicitly
%
%
\begin{eqnarray*}
I_d(\ell(t))
&=& \int_0^1 (1-[\mathbf e]_0)\log \frac{1-[\mathbf e]_0}{p(t)+(1-p(t))\frac{\beta(t)(e_0t+c_0)}{(1+\beta(t))t+\tilde c+c\beta(t)}}\\
&&\hspace*{13pt}{} + \sum_{i=1}^d (1-[\mathbf e]_i)\log \frac{1-[\mathbf e]_i}{(1-p(t))\frac{(i+\beta(t))(e_it+c_i)}{(1+\beta(t))t+\tilde c+c\beta(t)}}\\
&&\hspace*{13pt}{} + \Biggl(1-\sum_{i=0}^d (1-[\mathbf e]_i)\Biggr)\log
\frac{1-\sum_{i=0}^d (1-[\mathbf e]_i)}{(1-p(t))(1-\frac{\sum_{i=0}^d
(i+\beta(t))(e_it+c_i)}{(1+\beta(t))t+\tilde c+c\beta(t)})}\,dt
\end{eqnarray*}
is bounded under the bounds on $p,\beta$ in assumption (ND).
\end{pf}

\subsection{Lower bound}

Fix $h\dvtx C([0,1];\mathbb{R}^{d+2})\to\mathbb{R}$, a bounded, continuous
function, and $\varphi^*\in\Gamma_d$ such that $I_d(\varphi
^*)<\infty
$. To show the lower bound, it suffices to prove, for each $\varepsilon
>0$, that
%
%
\begin{equation}\label{eqLBgoal}
\limsup_{n\to\infty}V^n \leq I_d(\varphi^*)+h(\varphi^*)+8\eps.
\end{equation}
The main idea of the argument is to construct from $\varphi^*$ a
sequence of
control measures suitable to evaluate formulas for $V^n$.

Note only in this ``lower bound'' subsection, to make several
expressions simpler, we often take $c_{d+1}:= \bar c^d$.

\subsubsection{Step 1: Convex combination and regularization} Rather
than work
directly with $\varphi^*$, we consider a convex
combination of paths with better regularity: for $0\leq\theta\leq1$, let
\[
\varphi_\theta(t) = (1-\theta)\varphi^*(t) + \theta\ell(t),
\]
where $\ell(t)=\mathbf e t+\mathbf c^d$ is a linear function such that
$\mathbf e$ satisfies the assumptions of Lemma~\ref{Iforlinear}, say
$\mathbf{e}=(\frac12,\frac1{2^2},\ldots,\frac1{2^{d+1}},\frac1{2^{d+1}})$.
%
%
\begin{lemma}
As $\theta\downarrow0$, we have
\[
|I_d(\varphi_\theta) - I_d(\varphi^*)|\rightarrow0
\quad\mbox{and}\quad
|h(\varphi_\theta)- h(\varphi^*)| \rightarrow0.
\]
\end{lemma}
\begin{pf}
By convexity of $I_d$, and finiteness of $I_d(\ell(t))$ from Lemma
\ref
{Iforlinear},
\[
I_d(\varphi_\theta) \leq(1-\theta)I_d(\varphi^*) + \theta
I_d(\ell).
\]
On the other hand, since $|\varphi_\theta(t)-\varphi^*(t)| = |\int_0^t
(\dot\varphi_\theta-\dot\varphi^*)(s) \,ds| \leq2t\theta(d+2)$, we have
$\|\varphi_\theta-\varphi^*\|_\infty< 2\theta(d+2)\downarrow0$, by
lower semi-continuity of $I_d$, we have
\[
\liminf_{\theta\downarrow0} I_d(\varphi_\theta) \geq
I_d(\varphi^*).
\]

Also, as $h$ is continuous, we have
that $|h(\varphi_\theta)-h(\varphi^*)|\rightarrow0$.
\end{pf}

Now, fix $\theta>0$ such that
\[
I_d(\varphi_\theta) \leq I_d(\varphi^*) + \eps
\quad\mbox{and}\quad
h(\varphi_\theta) \leq h(\varphi^*) + \eps.
\]
Next,
for $\kappa\in\N$ and $t\in
[0,1]$, define
%
%
\begin{equation}\label{phikappa}
\psi_\kappa(t) = \int_0^t \gamma_\kappa(s) \,ds + \mathbf c^d,
\end{equation}
where
\[
\gamma_\kappa(t) = \kappa
\int_{i/\kappa}^{(i+1)/\kappa}\dot{\varphi_\theta}(s) \,ds %
\]
for $t\in[i/\kappa, (i+1)/\kappa)$, $0\leq i\leq\kappa-1$,
and $\gamma_\kappa(1) = \gamma_\kappa(1-1/\kappa)$. Note that
$\psi
_\kappa\in\Gamma_d$, and on
$[i/\kappa, (i+1)/\kappa)$ for
$0\leq i\leq\kappa-1$,
$\dot\psi_\kappa(t)$ equals the constant vector
$\gamma_\kappa(i/\kappa)$. In particular, $\dot\psi_\kappa$ is a
step function.
%
%
\begin{lemma}\label{propertiesofpsikappa}For $0\leq i\leq d+1$ and
$0\leq t\leq1$,
%
%
\begin{eqnarray}
\label{kappabounded}
\psi_{\kappa,i}(t) &\geq& \theta( e_i t+c_i),\\
\label{kappadotsumbounded}
\sum_{i=0}^{d+1} i\dot\psi_{\kappa,i}(t) &=& \sum_{i=0}^d
\bigl(1-[\dot\psi
_\kappa(t)]_i\bigr) \leq 1-\theta e_{d+1}.
\end{eqnarray}
\end{lemma}
\begin{pf}
These are properties of $\varphi_\theta$ inherited from properties of
$\varphi^*, \ell\in\Gamma_d$, which are preserved with respect to
(\ref{phikappa}).
Indeed, for each $0\leq i\leq d+1$,
\begin{eqnarray*}
\psi_{\kappa,i}(t)
&=& \varphi_{\theta,i}(\lfloor t\kappa\rfloor/\kappa) +(t\kappa
-\lfloor t\kappa\rfloor)\bigl(\varphi_{\theta,i}\bigl((\lfloor
t\kappa
\rfloor+1)/\kappa\bigr)
-\varphi_{\theta,i}(\lfloor t\kappa\rfloor/\kappa)\bigr)\\
&\geq& \theta(e_it+c_i).
\end{eqnarray*}
Last, (\ref{kappadotsumbounded}) follows: noting that
$\sum_{i=0}^d (1-[\mathbf e]_i)=\sum_{i=0}^{d+1} ie_i=1-e_{d+1}$,
\begin{eqnarray*}
&&\sum_{i=0}^d \bigl(1-[\dot\psi_\kappa(t)]_i\bigr)\\
&&\qquad= \kappa\int_{l/\kappa}^{(l+1)/\kappa} \Biggl[(1-\theta)\sum_{i=0}^d
\bigl(1-[\dot\varphi^*(s)]_i\bigr)
+ \theta\sum_{i=0}^d \bigl(1-[\dot\ell
(s)]_i\bigr)
\Biggr] \,ds\\
&&\qquad\leq 1-\theta+ \theta\sum_{i=0}^d (1-[\mathbf e]_i) = 1
-
\theta e_{d+1}.
\end{eqnarray*}
\upqed\end{pf}
%
%
\begin{lemma}
For large enough $\kappa$, we have
%
%
\begin{equation}\label{hIkappa}
h(\psi_\kappa) \leq h(\varphi^*)+2\eps
\quad\mbox{and}\quad I_d(\psi
_\kappa) \leq I_d(\varphi^*)+2\eps.
\end{equation}
\end{lemma}
\begin{pf}
Since
\[
\lim_{\kappa\rightarrow\infty}\sup_{t\in[0,1]}|\varphi_\theta
(t) -
\psi_\kappa(t)| = 0,
\]
the inequality with respect to $h$ follows from continuity of $h$ and
choosing $\kappa$ in terms of $\theta$.
We also note, by absolute continuity of $\varphi_\theta$, that a.s. in
$t$,
\[
\dot{\psi}_\kappa(t) = \gamma_\kappa(t) =
\kappa\int_{\lfloor t\kappa\rfloor/\kappa}^{(\lfloor t\kappa
\rfloor
+1)/\kappa} \dot{\varphi_\theta}(s) \,ds \rightarrow
\dot{\varphi_\theta}(t) \qquad\mbox{as
$\kappa\uparrow\infty$}.
\]

Then, by the form of $L$ [cf. (\ref{Lboundedcase})], bounds in Lemma
\ref{propertiesofpsikappa} and piecewise continuity and bounds on
$p,\beta$ in assumption (ND),\vspace*{1pt}
we have, as $\kappa\uparrow\infty$, that
$L(\vec\xi(t),\psi_\kappa(t),\dot\psi_\kappa(t))\to L(\vec\xi
(t),\varphi_\theta(t),\dot\varphi_\theta(t))$ for almost all
$t\in[0,1]$.

Also, we can bound $L(\vec\xi(t),\psi_\kappa(t),\dot\psi_\kappa(t))$
as follows: first, using $x\log x \leq0$ for $0\leq x\leq1$, bound that
\begin{eqnarray*}
&&L(\vec\xi(t),\psi_\kappa(t),\dot\psi_\kappa(t))\\
&&\qquad \leq- \bigl(1-[\dot{\psi_{\kappa}}(t)]_0\bigr)
\log\biggl(p(t)+\bigl(1-p(t)\bigr)\frac{\beta(t)\psi_{\kappa,0}(t)}{(1+\beta
(t))t+\tilde c+c\beta(t)}\biggr)\\
&&\qquad\quad{} - \sum_{i=1}^d \bigl(1-[\dot{\psi_{\kappa}}(t)]_i\bigr) \log
\biggl(\bigl(1-p(t)\bigr)\frac{(i+\beta(t))\psi_{\kappa,i}(t)}{(1+\beta
(t))t+\tilde c+c\beta(t)}\biggr)\\
&&\qquad\quad{} -\Biggl(1-\sum_{i=0}^d\bigl(1-
[\dot{\psi_{\kappa}}(t)]_i\bigr)\Biggr)\\
&&\qquad\quad\hspace*{11.2pt}{}\times
\log
\biggl(\bigl(1-p(t)\bigr)\biggl(1-\frac{\sum_{i=0}^d (i+\beta(t))\psi_{\kappa
,i}(t)}{(1+\beta(t))t+\tilde c+c\beta(t)}\biggr)\biggr).
\end{eqnarray*}
Now, as $0\leq[\dot\psi_\kappa]_i \leq1$ and $0\leq\sum_{i=0}^d(1-
[\dot{\psi_{\kappa}}]_i) \leq1$, we have the further upperbound,
using Lemma~\ref{propertiesofpsikappa},
\begin{eqnarray*}
&&-\log\biggl( {p(t)+\bigl(1-p(t)\bigr)\frac{\beta(t)\theta(e_0t+c_0)}{(1+\beta
(t))t+\tilde c+c\beta(t)}}\biggr)\\
&&\qquad{} - \sum_{i=1}^d\log\biggl({\bigl(1-p(t)\bigr)\frac{(i+\beta(t))\theta
(e_it+c_i)}{(1+\beta(t))t+\tilde c+c\beta(t)}}\biggr)\\
&&\qquad{} - \log\biggl( {\bigl(1-p(t)\bigr)\frac{(d+1+\beta(t))\theta
(e_{d+1}t+\bar c^d)}{(1+\beta(t))t+\tilde c+c\beta(t)}}\biggr),
\end{eqnarray*}
which is integrable on $[0,1]$ given the bounds on $p,\beta$ in
assumption (ND).

By dominated convergence, we obtain $\lim_\kappa I_d(\psi_\kappa)=
I_d(\varphi_\theta)$, and therefore the other inequality with respect
to $I_d$.
\end{pf}

Let now $\kappa$ be such that (\ref{hIkappa}) holds.
Finally,
we modify $\psi_\kappa$ on the interval $[0,\delta]$, for a small
enough $\delta>0$ to be chosen later.

Define
%
%
\begin{equation}\label{tieq}
t_i:= \delta-\sum_{l=i}^{d}\bigl(\delta+[\mathbf c^d]_l-[\psi
_\kappa
(\delta)]_l\bigr)
\end{equation}
for $0\leq i\leq d$, and $t_{d+1}:=\delta$; set also $t_{-1}:=0$ and
$t_{d+2}=t_{d+1}$.
Let also
%
%
\begin{equation}\label{psi^*}
\psi^*(t) =
\int_0^t \gamma^*(s) \,ds + \mathbf c^d,
\end{equation}
where
\[
\gamma^*(t)=
\cases{
\mathbf f_{d+1}, &\quad when $0 \leq t< t_{0}$,\cr
\mathbf f_i, &\quad when $t_{i} \leq t< t_{i+1}, 0\leq i\leq d$,\cr
\gamma_\kappa(t), &\quad when $t\geq\delta$.}
\]
Note $\gamma^*$ may not be defined at some endpoints as possibly $t_i =
t_{i+1}$ for some~$i$.

By inspection, $\psi^*\in\Gamma_d$.
Also,
$\dot\psi^*(t) = \mathbf f_{d+1}$ when $0\leq t<t_{0}$ and $\dot\psi
^*(t) = \mathbf f_i$ when $t_i\leq t<t_{i+1}$ for $0\leq i\leq d$.
Moreover, we have the following properties.

%
\begin{lemma}\label{psi^*forsmallt} We have
$\psi^*(\delta)=\psi_\kappa(\delta)$ and $t_0\geq\theta
e_{d+1}\delta$.
Also,
\[
\psi^*_0(t)=t+c_{0},\qquad
\psi^*_j(t)=c_j \qquad\mbox{for } 1\leq j\leq d+1,
\]
when $0\leq t<t_0$, and
\begin{eqnarray*}
\psi^*_0(t)&\geq&\theta e_{d+1}\delta+c_0\qquad \mbox{when }
t_0< t< t_1,\\
\psi^*_i(t)&\geq&\theta( e_i \delta+c_i)\qquad \mbox{when } t_i<
t< t_{i+1} \mbox{ and } 1\leq i\leq d.
\end{eqnarray*}
\end{lemma}
\begin{pf}
The lower bound for $t_0$ follows from the integration of both sides in
(\ref{kappadotsumbounded}) and the definition of $t_0$.
Now, we note that $\dot\psi^*_0(t)=0$ if $t_0\leq t\leq t_1$, and $1$
otherwise. Also, note that for $1\leq i\leq d+1$, $\dot\psi^*_i(t)=1$ if
$t_{i-1}<t< t_i$, $\dot\psi^*_i(t)=-1$ if $t_{i}< t< t_{i+1}$ and
$\dot\psi^*_i(t)=0$ otherwise.
Thus, noting (\ref{tieq}),
\[
\psi^*_0(\delta)=\int_0^\delta\gamma_0^*(s) \,ds + c_0 =
\delta- (t_1-t_0) + c_0=\psi_{\kappa,0}(\delta)
\]
and, for $1\leq i\leq d+1$,
\[
\psi^*_i(\delta) =\int_0^\delta\gamma_i^*(s) \,ds + c_i =
(t_i-t_{i-1}) - (t_{i+1}-t_i) + c_i = \psi_{\kappa,i}(\delta),
\]
which proves that $\psi^*(\delta)=\psi_\kappa(\delta)$.
Since $\psi^*_0(t)$ is nondecreasing, for $t\geq t_0$, $\psi
^*_0(t)\geq
\psi^*_0(t_0) = t_0+c_0\geq\theta e_{d+1}\delta+ c_0$.
For $1\leq i\leq d$, for $t_i< t<t_{i+1}$, $\psi^*_i(t)$ decreases to
its final value $\psi_{\kappa,i}(\delta)\geq\theta( e_i \delta+c_i)$
by (\ref{kappabounded}).
\end{pf}

\subsubsection{\texorpdfstring{Step 2: More properties of $\psi^*$}{Step 2: More properties of psi*}}
We now show the rate of $\psi^*$ up to time $\delta$ does not
contribute too much.
%
%
\begin{lemma}\label{handIforpsi^*}
For small enough $\delta>0$,
\[
\int^\delta_0 L(\vec\xi(t),\psi^*(t),\dot
\psi^*(t)) \,dt \leq\varepsilon
\quad\mbox{and}\quad \|\psi^* -
\psi
_\kappa\|_\infty<\varepsilon.
\]
In particular,
$h(\psi^*) \leq h(\varphi^*) + 3\varepsilon$ and $I_d(\psi^*) \leq
I_d(\varphi^*) + 3\varepsilon$.
\end{lemma}
\begin{pf}
\hspace*{-0.5pt}Write, for $0\leq t\leq\delta$, as $L(\vec\xi(t),\hspace*{-0.5pt}\psi^*(t),
\hspace*{-0.5pt}\dot\psi^*(t))
= R(\delta_{\mathbf f_{d+1}} \hspace*{-0.5pt}\mmid\rho_{\vec\xi(t),\psi
^*(t)})
\hspace*{-0.5pt}\times1(0<t<t_0) + \sum_{i=0}^{d} R(\delta_{\mathbf
f_{i}} \mmid\rho_{\vec\xi(t),\psi^*(t)})1(t_i<t<t_{i+1})$,
\begin{eqnarray*}
&&L(\vec\xi(t),\psi^*(t),\dot\psi^*(t))\\
&&\qquad =
- \log\biggl(\bigl(1-p(t)\bigr)\biggl(1-\frac{\sum_{l=0}^d(l+\beta(t))\psi
^*_l(t)}{(1+\beta(t))t+\tilde c+c\beta(t)}\biggr)\biggr)1(0<t<t_0)\\
&&\qquad\quad{} - \log\biggl(p(t) + \bigl(1-p(t)\bigr)\frac{\beta(t)\psi
^*_0(t)}{(1+\beta(t))t+\tilde c+c\beta(t)}\biggr)1(t_0<t<t_{1})\\
&&\qquad\quad{} - \sum_{i=1}^d\log\biggl(\bigl(1-p(t)\bigr)\frac{(i+\beta
(t))\psi
^*_i(t)}{(1+\beta(t))t+\tilde c+c\beta(t)}\biggr)1(t_i<t<t_{i+1}).
\end{eqnarray*}
By Lemma~\ref{psi^*forsmallt} and the bounds on $p,\beta$ in
assumption (ND), this expression is integrable for $0\leq t\leq\delta$.
(It would be bounded unless $\bar c^d=0$ and $c\neq0$, in which case
the first term in the expression involves $-\log t$.)
Hence, the first statement follows for small $\delta>0$.
Also, the second statement holds as $\|\psi^*-\psi_\kappa\|_\infty=
\sup_{0\leq t<\delta} |\psi^*-\psi_\kappa|\leq2\delta(d+2)$.
The last\vspace*{1pt} statement is a consequence now of (\ref{hIkappa}).
\end{pf}

We will take $\delta>0$ small enough so that the bounds in the above
lemma hold.
%
%
\begin{lemma}
We have
%
%
\begin{equation}\label{psisup}
\lim_{n\to\infty}\sup_{0\leq j\leq n}\Biggl|\psi^*(j/n)-\frac1 n
\sum_{l=0}^{j-1}\dot\psi^*(l/n)-\mathbf c^d\Biggr| = 0.
\end{equation}
Also, for
$j\geq\lfloor\delta n\rfloor$ and $0\leq i\leq d+1$,
%
%
\begin{equation}\label{LBpsi}
\frac{1}{n}\sum_{l=0}^{j-1} \dot{\psi}^*_{i} (l/n)+c_i
\geq\frac\theta2\biggl(\frac{e_i j}{n}+c_i\biggr).
\end{equation}
\end{lemma}
\begin{pf} Since $\dot\psi^*$ is piecewise constant, when $l/n\leq
s\leq(l+1)/n$,
$|\dot\psi^*(s)-\dot\psi^*(l/n)|\ne0$ for at most
$\kappa$\vadjust{\goodbreak}
subintervals [cf. (\ref{phikappa})
and (\ref{psi^*})], and is also bounded by
$2(d+2)$. Hence,
\begin{eqnarray*}
\Biggl|\psi^*(j/n)-\frac1 n \sum_{l=0}^{j-1}\dot\psi^*(l/n)-\mathbf
c^d\Biggr|
&=&\Biggl|\sum_{l=0}^{j-1}\int_{l/n}^{(l+1)/n}\bigl(\dot\psi
^*(s)-\dot\psi
^*(l/n)\bigr) \,ds\Biggr|\\[-2pt]
&\leq& \frac{2(d+2)} n \kappa.
\end{eqnarray*}
The last statement follows from (\ref{kappabounded}).
\end{pf}

\subsubsection{Step 3: Admissible control measures and convergence} We
now build a sequence of
controls based on $\psi^*$.
Define $\nu_0=\nu_0(\dot\psi^*(j/n)|j/n)$ using (\ref{nuforLbounded}), and
\begin{eqnarray*}
&&v^n_j(d\mathbf y;\mathbf{x}_{0},\ldots,\mathbf{x}_j)\\[-2pt]
&&\qquad =
\cases{
\nu_0\bigl(\dot\psi^*(j/n)|j/n\bigr),
&\quad when $0\leq j\leq\lfloor\delta n\rfloor$\vspace*{2pt}\cr
&\qquad or when $j\geq\lceil\delta n\rceil$\vspace*{2pt}\cr
&\qquad and $\mathbf{x}_{j,i}\geq\dfrac\theta
4(e_i\delta+c_i)$\vspace*{2pt}\cr
&\qquad for $0\leq i\leq d+1$,\cr
\rho_{\vec\xi(j/n), \mathbf{x}_j}, &\quad otherwise.}
\end{eqnarray*}
The reasoning behind this choice of controls is as follows: to bound
the limit of the quantity in (\ref{eqrepformulafornonleaves}), using
formula (\ref{Lboundedcase}), by $I_d(\psi^*) + h(\psi^*)$, we would
like to specify the controls in form $\nu_0(\dot\psi^*(j/n)|j/n)$.
Such a\vspace*{1pt} choice, as we will see, also ensures that the adapted sequence
$\bar{\mathbf X}^n(j)$ is close to $\psi^*(j/n)$. However, the adapted
process, as it is random, may get too close to a boundary. When this
happens, not often it turns out, to bound errors, we specify that the
controls take the cost-free form of the natural evolution sequence.
Also, to get past this boundary layer initially, $\psi^*$ has been
built as a step function so that the adapted process must follow a
deterministic trajectory up to time~$\lfloor\delta
n\rfloor$.\vspace*{1pt}\looseness=-1

Define $\bar{\mathbf X}^{n}(0) =\mathbf c^d$, and $\bar{\mathbf
X}^{n}(j+1) = {\bar{\mathbf X}^{n}(j)} +\frac1n{\bar{\mathbf
Y}^{n}(j)}$ for $j\geq0$ where
\[
\bar{P}\bigl({\bar{\mathbf Y}^{n}(j)}\in d\mathbf{y}|\bar{\mathbf
X}^{n}(0),\ldots, {\bar{\mathbf X}^{n}(j)}\bigr) =
v_j^n(d\mathbf{y};\bar{\mathbf X}^{n}(0),\ldots,
{\bar{\mathbf X}^{n}(j)}).
\]
Thus, for $j\geq0$, $\bar{\mathbf X}^{n}(j)=\frac
1n\sum_{l=0}^{j-1}\bar{\mathbf Y}^{n}(l)+\mathbf c^d$. It will be
useful later to note the
total weight $\sum_{i=0}^{d+1}
(i+\beta(j/n))\bar{\mathbf X}^{n}_{i}(j)\leq(j/n +\tilde c)
+\beta(j/n)(j/n +c)$ and, for $0\leq j\leq\lfloor\delta n\rfloor$, as
mentioned $\bar{\mathbf X}^{n}(j)$ is
deterministic and $\bar{\mathbf X}^{n}(j)=\break\frac1n\sum
_{l=0}^{j-1}\dot
\psi^*(l/n)+\mathbf c^d$.

Define now, for each $n\geq1$, the martingale sequence for $0\leq
j\leq n$
\begin{eqnarray*}
\mathbf M^{n}(j) :\!&=& \frac1 n \sum_{l=0}^{j-1} \bigl(\bar{\mathbf
Y}^{n}(l)-\bar E(\bar{\mathbf Y}^{n}(l)|\bar{\mathbf
X}^{n}(l))\bigr)\\[-2pt]
& = & \bar{\mathbf X}^{n}(j) - \frac1 n \sum_{l=0}^{j-1}
\bar E(\bar{\mathbf Y}^{n}(l) |\bar{\mathbf X}^{n}(l))
-
\mathbf c^d.
\end{eqnarray*}
Let
\begin{eqnarray*}
\tau_n &:=& n \wedge\min\biggl\{\lceil\delta n\rceil\leq l \leq n\dvtx
\bar{\mathbf X}^{n}_{i}(l)< \frac\theta4(e_{i}\delta+ c_{i})\\
&&\hspace*{91pt}
\mbox{for some $0\leq i\leq d+1$}\biggr\}.%
\end{eqnarray*}
Then, $\tau_n\geq\lceil\delta n\rceil$ is a stopping time, and the
corresponding stopped
process $\{\mathbf M^n(j\wedge\tau_n)\}$ is also a martingale for
$0\leq j\leq n$.
Let now
\[
\mathbb A_n:=
\biggl\{\sup_{0\leq j\leq n}|\mathbf M^n(j\wedge\tau_n)|>
\frac{\theta e_{d+1}}{4n^{1/8}}\biggr\}.
\]

%
\begin{lemma}
\label{taun=n}
For $n\geq\delta^{-8}$, on the set $\mathbb A_n^{c}$, we have $\tau_n
= n$.
\end{lemma}
\begin{pf} From the definition of $\{v_j^n\}$ and $\tau_n$, we have
$\bar E(\bar{\mathbf Y}^n(l)|\bar{\mathbf X}^n(l)) = \dot\psi^*(l/n)$
for $0\leq l\leq j\wedge\tau_n -1$ and $j\geq\lceil\delta n\rceil$.
Then, on $\mathbb A_n^{c}$, by (\ref{LBpsi}), we have
\begin{eqnarray*}
\bar{\mathbf X}^{n}_i(j\wedge\tau_n) &\geq& c_i +
\frac{1}{n}\sum_{l=0}^{j\wedge\tau_{n}-1} \bar E(\bar{\mathbf
Y}^{n}_{i}(l) |\bar{\mathbf X}^{n}(l))- \frac{\theta
e_{d+1}}{4n^{1/8}}
\\
&=& c_i+ \frac1n\sum_{l=0}^{j\wedge\tau
_n-1}\dot\psi^*_i(l/n) - \frac{\theta e_{d+1}}{4n^{1/8}}\\
&\geq&\frac\theta2\biggl(\frac{e_i(j\wedge
\tau
_n)}{n}+c_i\biggr) -\frac{\theta e_{d+1}}{4n^{1/8}}\\
&\geq&\frac\theta4(e_i\delta+ c_i).
\end{eqnarray*}
Hence, $\tau_n = n$.
\end{pf}

We now observe, by Doob's martingale inequality and bounds, in terms of
constants $C=C_d$, that
%
%
\begin{eqnarray}\label{Rosenthal}
\bar{P}[\mathbb A_n]
&\leq& Cn^{1/2}\bar{E}|\mathbf M^n(j\wedge\tau_n)
|^4\nonumber\\
&=& Cn^{-7/2}\bar{E}\Biggl|\sum_{l=0}^{j\wedge\tau_n-1}
\bigl(\bar{\mathbf Y}^{n}(l)
-\bar E(\bar{\mathbf Y}^{n}(l)|\bar{\mathbf X}^{n}(l))\bigr)\Biggr|^4\\
& \leq& Cn^{-7/2} n^2 = Cn^{-3/2}.\nonumber
\end{eqnarray}
We now state the following almost sure convergence.
%
%
\begin{lemma} We have
%
%
\begin{equation}\label{as}
\lim_{n\uparrow\infty} \sup_{0\leq j\leq n}
\Biggl|{\bar{\mathbf X}^n(j)}-\frac1 n \sum_{l=0}^{j-1}\dot{\psi
^*}(l/n)-\mathbf c^d\Biggr| = 0\qquad
\mbox{a.s.}
\end{equation}
\end{lemma}
\begin{pf}
First, by (\ref{Rosenthal}) and the Borel--Cantelli lemma, $\bar
{P}(\limsup\mathbb A_n) = 0$.
On the other hand, on the full measure set $\bigcup_{m\geq1} \bigcap
_{k\geq
m} \mathbb A_k^c$, since $\tau_n =n$ and $\bar E(\bar{\mathbf
Y}^{n}(l) |\bar{\mathbf X}^{n}(l))=\dot\psi^*(l/n)$ for
$0\leq
l\leq n-1$ on $\mathbb A_n^c$ by Lemma~\ref{taun=n}, the desired
convergence holds.%
\end{pf}

\subsubsection{Step 4} We now argue the lower bound through
representation (\ref{eqrepformulafornonleaves}). Recall the
definition of $\vec\xi(\cdot)$ in the beginning of Section \ref
{sectionLDP}.
The sum in
(\ref{eqrepformulafornonleaves}) equals
%
%
\begin{eqnarray}\label{splitUB}
&&\bar{E}\Biggl[\frac{1}{n} \sum_{j=0}^{n-1}
R\bigl(v^n_j \mmid\rho_{\vec\xi(j/n), \bar{\mathbf
X}^{n}(j)}\bigr)\Biggr]\nonumber\\
&&\qquad = \bar{E}\Biggl[ \frac{1}{n} \sum_{j=0}^{\lfloor\delta
n\rfloor} R\bigl(v^n_j \mmid\rho_{\vec\xi(j/n), \bar{\mathbf
X}^{n}(j)}\bigr)
\Biggr]\nonumber\\
&&\qquad\quad{} + \bar{E}\Biggl[ \frac{1}{n} \sum^{n-1}_{j=\lceil\delta
n\rceil} R\bigl(v^n_j \mmid\rho_{\vec\xi(j/n), \bar{\mathbf X}^{n}(j)}\bigr);
\mathbb A_n \Biggr]\\
&&\qquad\quad{} + \bar{E}\Biggl[ \frac{1}{n} \sum_{j=\lceil\delta
n\rceil}^{n-1} R\bigl(v^n_j \mmid\rho_{\vec\xi(j/n), \bar{\mathbf X}^{n}(j)}\bigr);
\mathbb A_n^c\Biggr]\nonumber\\
&&\qquad = A_1 + A_2 + A_3.\nonumber
\end{eqnarray}

\textit{Step} 4.1. We treat the term $A_2$ in (\ref{splitUB}).
Recall $\sigma(j/n)=(1+\beta(j/n))(j/n) + \tilde c + c\beta(j/n)$ and
the ``weight'' bound on $\bar{\mathbf X}^n(j)$ in beginning of Step 3.
For $\lceil\delta n\rceil\leq j\leq n-1$,
\begin{eqnarray*}
&&
R\bigl(v^n_j \mmid\rho_{\vec\xi(j/n), \bar{\mathbf X}^{n}(j)}\bigr)\\
&&\qquad=
R\bigl(\nu_0\bigl(\dot{\psi^*}(j/n)\bigr) \mmid\rho_{\vec\xi(j/n), \bar
{\mathbf
X}^{n}(j)}\bigr)\\
&&\qquad\quad{} \times1\bigl(\bar{\mathbf X}^{n}_{i}(j)\geq(\theta
/4)(e_i\delta+c_i) \mbox{ for }0\leq i\leq d+1\bigr).
\end{eqnarray*}
Noting (\ref{Lboundedcase}), this is bounded above, using $x\log x
\leq0$ for $0\leq x\leq1$, by
\begin{eqnarray*}
&&\Biggl[- \biggl(1-\biggl[\dot\psi^*\biggl(\frac jn\biggr)\biggr]_0\biggr)
\log\biggl(p(j/n) + \bigl(1-p(j/n)\bigr)\frac{\beta(j/n)\bar{\mathbf
X}^{n}_0(j)}{\sigma(j/n)}\biggr)\nonumber\\
&&\qquad{} - \sum_{i=1}^d\biggl(1-\biggl[\dot\psi^*\biggl(\frac jn\biggr)
\biggr]_i\biggr)
\log\biggl(\bigl(1-p(j/n)\bigr)\frac{(i+\beta(j/n))\bar{\mathbf
X}^{n}_i(j)}{\sigma(j/n)}\biggr)\nonumber\\
&&\qquad{} - \Biggl(\sum_{i=0}^d\biggl[\dot\psi^*\biggl(\frac jn\biggr)\biggr]_i
- d\Biggr)\\
&&\hspace*{63.5pt}{}\times
\log\biggl(\bigl(1-p(j/n)\bigr)\biggl(1-\frac{\sum_{i=0}^d(i+\beta
(j/n))\bar
{\mathbf X}^{n}_{i}(j)}
{\sigma(j/n)}\biggr)\biggr)\Biggr]\\
&&\hspace*{11pt}\qquad{} \times1\bigl(\bar{\mathbf
X}^n_{i}(j)\geq
(\theta/4)(e_i\delta+c_i) \mbox{ for }0\leq i\leq d+1
\bigr).%
\end{eqnarray*}
Given bounds on $p,\beta$ in (ND), as $0\leq[\dot\psi^*]_i \leq1$, we
have $d\leq\sum_{i=1}^d[\dot\psi^*]_i \leq d+1$ and
\begin{eqnarray*}
\sum_{i=0}^d \bigl(i+\beta(j/n)\bigr)\bar{\mathbf X}^n_{i}(j)
&\leq&\sigma(j/n) - \bigl(d+1+\beta(j/n)\bigr)\bar{\mathbf
X}^n_{d+1}(j)\\
&\leq&\sigma(j/n) - \bigl(d+1+\beta(j/n)\bigr)\cdot
(\theta
/4)(e_{d+1}\delta+ c_{d+1}),
\end{eqnarray*}
the relative entropy is further bounded by a constant $C_d$.
Thus, for large $n$,
%
%
\begin{equation}\label{A2}
A_2
\leq C_d\cdot
\bar{P}\biggl[\sup_{0\leq j\leq n}|\mathbf M^n(j\wedge\tau_n)|>
\frac{\theta e_{d+1}}{4n^{1/8}}\biggr] \leq\eps.
\end{equation}

\textit{Step} 4.2. Now, for the term $A_1$ in (\ref{splitUB}),
we recall for $j\leq\lfloor\delta n\rfloor$ that
$\bar{\mathbf X}^n(j)=\frac1n\sum_{l=0}^{j-1}\dot\psi
^*(l/n)+\mathbf
c^d$ is deterministic.
Also note, for $0\leq i\leq d$, that $\dot{\psi^*}(t)=\mathbf f_i$ on
$t_{i}< t< t_{i+1}$, and $\dot{\psi^*}(t)=\mathbf f_{d+1}$ on
$0=t_{-1}\leq t\leq t_0$ (cf. near Lemma~\ref{psi^*forsmallt}).
Thus, for $0\leq j\leq\lfloor\delta n\rfloor$,
denoting $\mathbf f_{-1}=\mathbf f_{d+1}$, we may write
\begin{eqnarray*}
&&
R\bigl(v^n_j \mmid\rho_{\vec\xi(j/n), \bar{\mathbf X}^n(j)}\bigr) \\
&&\qquad= L
\Biggl(\vec
\xi\biggl(\frac jn\biggr),\frac{1}{n}\sum_{l=0}^{j-1}\dot\psi^*
\biggl(\frac
ln\biggr)+\mathbf c^d, \dot{\psi^*}\biggl(\frac jn\biggr)\Biggr)\\
&&\qquad= \sum_{i=-1}^d L\Biggl(\vec\xi\biggl(\frac jn\biggr),\sum
_{l=-1}^{i-1}\frac{\lfloor t_{l+1} n\rfloor- \lfloor t_{l} n\rfloor
}{n}\mathbf f_l + \frac{j-\lfloor t_{i} n\rfloor}{n}\mathbf f_i +
\mathbf c^d, \mathbf f_i\Biggr)\\
&&\hspace*{19.5pt}\qquad\quad{}
\times1(\lfloor t_i n\rfloor\leq j<\lfloor t_{i+1} n\rfloor),
\end{eqnarray*}
where, for $i=-1$, the empty sum in the argument for $L$ vanishes.
Comparing with the proof of Lemma~\ref{handIforpsi^*}, this
expression, given bounds on $p,\beta$ in (ND), is
bounded, for $0\leq j\leq\lfloor\delta n\rfloor$, except when $\bar
c^d=0$ and $c\neq0$, in which case a
``$-\log(j/n)$'' term appears in the $i=-1$ term. But, since
$-(1/n)\sum_{j=1}^{\lfloor\delta n\rfloor} \log(j/n) \leq
-\int_0^\delta\log(t)\,dt$, its contribution is still small. Hence,
%
%
\begin{equation}\label{A1}
A_1
\leq\epsilon(\delta)\qquad \mbox{where }\epsilon(\delta)\rightarrow0 \mbox{ as } \delta\to0.
\end{equation}
%

\textit{Step} 4.3. We now estimate the last term $A_3$ in (\ref
{splitUB}). For $n\geq\delta^{-8}$, by Lem\-ma~\ref{taun=n} and
definition of $L$ (\ref{Lboundedcase}),
\begin{eqnarray*}
A_3
&\leq& \bar{E}\Biggl[\frac{1}{n}\sum_{j=\lceil\delta n\rceil}^{n-1}
L\biggl(\vec\xi\biggl(\frac jn\biggr),\bar{\mathbf X}^n(j), \dot{\psi
^*}
\biggl(\frac jn\biggr)\biggr); \mathbb A_n^c \cap\{\tau_n=n\} \Biggr]\\
&\leq& \bar{E}\biggl[\int_{\delta}^{1} L\biggl(\vec\xi\biggl(\frac
{\lfloor
nt\rfloor}n\biggr),\bar{\mathbf X}^n(\lfloor nt\rfloor), \dot{\psi
^*}
\biggl(\frac{\lfloor nt\rfloor}n\biggr)\biggr)\,dt; \mathbb A_n^c\cap\mathbb
B_n\biggr],
\end{eqnarray*}
where $\mathbb B_n = \{\bar{\mathbf X}_{i}^n(j)\geq(\theta
/4)(e_i\delta
+c_i) \mbox{ for } 0\leq i\leq d+1, j\geq\lceil\delta n\rceil\}$.
On the event $\mathbb A_n^c\cap\mathbb B_n$, $\bar E(\bar{\mathbf
Y}^n(l)|\bar{\mathbf X}^n(l)) = \dot\psi^*(l/n)$ for $l\geq0$, and so
$|\bar{\mathbf X}^n(\lfloor nt\rfloor) - \psi^*(\lfloor
nt\rfloor/n)|1(\mathbb A_n^c\cap\mathbb B_n)\rightarrow0$
for each realization from (\ref{psisup}) and (\ref{as}).

Also, from the form of $L$ (\ref{Lboundedcase}), as $\dot\psi^*$
is a
step function, (\ref{kappabounded}), and bounds and piecewise
continuity of $p,\beta$ in (ND), we may bound as in Step 4.1 and observe
\begin{eqnarray*}
2C_d &\geq& \biggl| L\biggl(\vec\xi\biggl(\frac{\lfloor
nt\rfloor}n\biggr),\bar{\mathbf X}^n(\lfloor nt\rfloor),
\dot{\psi^*}\biggl(\frac{\lfloor nt\rfloor}n\biggr)\biggr)\\
&&\hspace*{56.5pt}{}
- L(\vec\xi(t),\psi^*(t),\dot\psi^*(t))\biggr|
1(\mathbb A_n^c\cap\mathbb
B_n) \\
&\rightarrow&0
\end{eqnarray*}
for almost all $t$ and each realization.
Hence, by bounded convergence theorem, with respect to
$d\bar{P}\times1([\delta,t])\,dt$,
%
%
\begin{equation}\label{A3}
\limsup_{n\rightarrow\infty} A_3
\leq\int_\delta^1 L(\vec\xi(t),\psi^*(t),\dot\psi
^*(t))\,dt.
\end{equation}

\subsubsection{Step 5} Finally,\vspace*{1pt} by (\ref{psisup}) and (\ref{as}),
$\lim_{n\rightarrow\infty}h(\bar{\mathbf X}^n(\cdot)) =
h(\psi^*(\cdot))$ a.s. in the sup topology, and by bounded convergence
$\lim_{n\rightarrow
\infty}\bar{E}[h(\bar{\mathbf X}^n(\cdot))] =
h(\psi^*(\cdot))$.

We now combine all bounds to conclude the proof of (\ref{eqLBgoal}).
By (\ref{eqrepformulafornonleaves}),
bounds (\ref{A2}), (\ref{A1}), (\ref{A3}) and nonnegativity of
$L$, we have
\begin{eqnarray*}
\limsup_{n\to\infty} V^n&\leq&\limsup_{n\to\infty}
\bar{E}\Biggl[\frac{1}{n}\sum_{j=0}^{n-1}
R\bigl(v^n_j \mmid\rho_{\vec\xi(j/n), \bar{\mathbf X}^n(j)}\bigr) + h(\bar
{\mathbf X}^n(\cdot))\Biggr]\\
&\leq& 2\eps+\int_0^1 L(\vec\xi(t),\psi^*(t),\dot\psi
^*(t))
\,dt + h(\psi^*).
\end{eqnarray*}
Then, by Lemma~\ref{handIforpsi^*}, we obtain
(\ref{eqLBgoal}).


\section{\texorpdfstring{Proof of Theorem \protect\ref{LDPinfin}}{Proof of Theorem 1.4}}
\label{sectionLDPinfin}

The proof of Theorem~\ref{LDPinfin} follows from the following two
propositions, and is given below. We first recall the projective limit
approach, following notation in~\cite{Dembo-Zeitouni}, Section 4.6.
Define, for $0\leq i\leq j$, $\mathcal Y_j = C([0,1]; \R^{j+2})$ and
$p_{ij}\dvtx\mathcal Y_j \rightarrow\mathcal Y_i$ by $\langle\varphi
_0,\ldots,
\varphi_{j+1}\rangle\mapsto\langle\varphi_0,\ldots, \varphi_i,\break
\sum_{l={i+1}}^{j+1} \varphi_l\rangle$. Also define $\varprojlim
\mathcal
Y_j \subset\prod_{i\geq0}\mathcal Y_i$ as the subset of elements $x=
\langle x^0,\break x^1,\ldots\rangle$ such that $p_{ij}x^j = x^i$, equipped
with the
product topology. Let also $p_j\dvtx\varprojlim\mathcal Y_j \rightarrow
\mathcal Y_j$ be the canonical projection, $p_j x = x^j$.

Since $I_d$ are convex, good rate functions on $C([0,1],\R^{d+2})$, by
the LDPs Theorem~\ref{LDPmain} and~\cite{Dembo-Zeitouni}, Theorem
4.6.1, we obtain the following proposition. Recall the
notation in Theorem~\ref{LDPmain}. For $n\geq1$, let $\mathcal
{X}^{n,\infty} = \langle\mathbf X^{n,0},\break \mathbf
X^{n,1},\ldots\rangle$.

%
\begin{prop}
\label{projectivelimit}
The sequence $\{\mathcal{X}^{n,\infty}\}\subset\varprojlim\mathcal
Y_j$ satisfies an LDP with rate $n$ and convex, good rate function
\[
J^\infty(\varphi) = \cases{
\displaystyle \sup_d \{I_d(p_d(\varphi))\}, &\quad when
$\varphi
\in\varprojlim\mathcal Y_j$,\vspace*{2pt}\cr
\infty, &\quad otherwise.}
\]
\end{prop}

To establish Theorem~\ref{LDPinfin}, it remains to further identify
$J^\infty$. Recall $\Gamma_d\subset C([0,1];\R^{d+2})$ are those
elements $\varphi= \langle\varphi_0,\ldots, \varphi_d,\varphi
_{d+1}\rangle$
such that:
\begin{itemize}
\item[]
$\varphi(0)=\mathbf c^d$, each\vspace*{1pt} $\varphi_i\geq0$ is Lipschitz with
constant 1 such that $0\leq[{\dot\varphi}(t)]_i\leq1$ for $0\leq
i\leq d$, $\sum_{i=0}^{d+1} \dot\varphi_i(t)=1$, and $\sum_{i=0}^{d+1}
i \dot\varphi_i(t)= \sum_{i=0}^d (1-[\dot\varphi(t)]_i)\leq1$ for
almost all $t$.
\end{itemize}
Let also $\Gamma^*\subset\varprojlim\mathcal Y_j$ be those elements
$\varphi= \langle\varphi^0,\varphi^1,\ldots\rangle$ such that $\varphi
^d \in
\Gamma_d$ for $d\geq0$. Since $\{\Gamma_d\}_{d\geq0}$ are compact
sets, it is a straightforward exercise to see that $\Gamma^*$ is compact.
Define $L_d(p_d(\varphi(t)))$ equal to
%
%
\begin{eqnarray*}
&&\bigl(1-[\dot{\varphi}^d(t)]_0\bigr) \log \frac{1-[\dot{\varphi}^d(t)]_0}{p(t)+(1-p(t))\frac{\beta(t)\varphi^d_0(t)}{(1+\beta(t))t+\tilde c+c\beta(t)}}\\
&&\qquad{} + \sum_{i=1}^d \bigl(1-[\dot{\varphi}^d(t)]_i\bigr) \log \frac{1-[\dot{\varphi}^d(t)]_i}{(1-p(t))\frac{(i+\beta(t))\varphi^d_i(t)}{(1+\beta(t))t+\tilde c+c\beta(t)}}\\
&&\qquad{} + \Biggl(1-\sum_{i=0}^d\bigl(1- [\dot{\varphi}^d(t)]_i\bigr)\Biggr) \log
\frac{1-\sum_{i=0}^d(1- [\dot{\varphi}^d(t)]_i)}{(1-p(t))(1-\frac{\sum_{i=0}^d (i+\beta(t))\varphi^d_i(t)}{(1+\beta(t))t+\tilde c+c\beta(t)})}\,dt.
\end{eqnarray*}

%
\begin{prop}
\label{rateinfinbounded}
The rate function $J^\infty(\varphi)$ diverges when $\varphi\notin
\Gamma^*$. However, for $\varphi\in\Gamma^*$,
$\lim_{d\uparrow\infty} L_d(p_d(\varphi(t)))$ exists for
almost all $t$, and we can evaluate
\[
J^\infty(\varphi) = \int_0^1 \lim_{d\uparrow\infty} L_d
(p_d(\varphi(t))) \,dt.
\]
\end{prop}
\begin{pf}
First, from the definition, $J^\infty(\varphi)$ diverges unless
$\varphi
\in
\Gamma^*$.
Next, for $\varphi\in\Gamma^*$ and almost all $t$,
we argue
%
%
\begin{equation}\label{rs} L_r(p_r(\varphi(t))) \leq
L_{s}(p_{s}(\varphi(t)))\qquad
\mbox{when }r<s.
\end{equation}
It will be enough to show
from the form of the rates the following:
%
%
\begin{eqnarray*}
&&\Biggl(1-\sum_{i=0}^r\bigl(1-[\dot\varphi^s(t)]_i\bigr)\Biggr)\log \frac{1-\sum_{i=0}^r(1-[\dot\varphi^s(t)]_i)}{(1-p(t))(1-\frac{\sum_{i=0}^r(i+\beta(t))\varphi^s_i(t)}{(1+\beta(t))t + \tilde c + c\beta(t)})} \\
&&\qquad\leq \sum_{i=r+1}^s \bigl(1-[\dot\varphi^s(t)]_i\bigr)\log \frac{1-[\dot\varphi^s(t)]_i}{(1-p(t))\frac{(i+\beta(t))\varphi^s_i(t)}{(1+\beta(t))t + \tilde c + c\beta(t)}} \\
&&\qquad\quad{} + \Biggl(1-\sum_{i=0}^s(1-[\dot\varphi^s(t)]_i)\Biggr)\log
\frac{1-\sum_{i=0}^s(1-[\dot\varphi^s(t)]_i)}{(1-p(t))(1-\frac{\sum_{i=0}^s(i+\beta(t))\varphi^s_i(t)}{(1+\beta(t))t
+ \tilde c + c\beta(t)})}.
\end{eqnarray*}

Consider now $h(x)= x\log x$ which is convex for $x\geq0$. Under
conventions (\ref{convention}), for nonnegative
numbers, $a_i$ and $b_i$, we have
\begin{eqnarray*}
\frac{\sum_{i=p}^qa_i}{\sum_{i=p}^qb_i}\log\frac{\sum
_{i=p}^qa_i}{\sum
_{i=p}^q b_i} &=& h\biggl(\frac{\sum_{i=p}^q a_i}{\sum
_{i=p}^qb_i}
\biggr) = h\Biggl(\sum_{i=p}^q \frac{b_i}{\sum_{i=p}^q b_i}\frac
{a_i}{b_i}\Biggr)\\
&\leq& \sum_{i=p}^q \frac{b_i}{\sum_{i=p}^q b_i}h\biggl(\frac
{a_i}{b_i}\biggr) = \frac{\sum_{i=p}^q a_i\log(a_i/b_i)}{\sum
_{i=p}^q b_i}.
\end{eqnarray*}
We now finish the proof of (\ref{rs}) by applying the last sequence,
with $p=r+1$ and $q=s+1$, to
\[
a_j = \cases{
\displaystyle 1-[\dot\varphi^s(t)]_j, &\quad for $r+1\leq
j\leq s$,\vspace*{2pt}\cr
\displaystyle 1-\sum_{i=0}^s \bigl(1-[\dot\varphi^s(t)]_i\bigr), &\quad for $j=s+1$,}
\]
and
\[
b_j = \cases{
\displaystyle \bigl(1-p(t)\bigr)\frac{(j+\beta(t))\varphi^s_j(t)}{(1+\beta
(t))t + \tilde c + c\beta(t)}, &\quad for $r+1\leq j\leq s$,\vspace*{2pt}\cr
\displaystyle \bigl(1-p(t)\bigr)\biggl(1-\frac{\sum_{i=0}^s(i+\beta(t))\varphi
^s_i(t)}{(1+\beta
(t))t + \tilde c + c\beta(t)}\biggr), &\quad for $j=s+1$.}
\]

Finally, given $L_d(p_d(\varphi(t)))\geq0$ is increasing in $d$, the
identification of $J^\infty$ in the display of the proposition follows
from monotone convergence.
\end{pf}
\begin{pf*}{Proof of Theorem~\ref{LDPinfin}}
Let $\Gamma^{\infty}\subset\prod_{i\geq0}C([0,1]; \R)$, endowed with
the product topology, be those elements $\xi= \langle\xi_0,\xi_1,\ldots
\rangle$
such that:
\begin{itemize}
\item[] $\xi_i(0)=c_i$, $\xi_i(t)\geq0$ is\vspace*{1pt} Lipschitz with constant
$1$, $0\leq[{\dot\xi}(t)]_i\leq1$ for $i\geq0$, and
$\frac{d}{dt}\sum_{i\geq0}\xi_i(t) = 1$ and\vspace*{1pt} $\lim_d [\sum_{i=0}^d
i\dot\xi_i(t) + (d+1)(1-[\dot\xi(t)]_d)]=\sum_{i\geq0} (1-[\dot
\xi
(t)]_i) \leq1$ for almost all $t$.
\end{itemize}

We now show that $\Gamma^\infty$ and $\Gamma^*$ are homeomorphic.
Hence, as $\Gamma^*$ is compact, $\Gamma^\infty$ would also be compact.
(We note, one can see directly that $\Gamma^\infty$ is compact.)

Define the map $F\dvtx\Gamma^\infty\rightarrow\Gamma^*$ by
\[
F(\xi) = \langle\xi^0,\ldots, \xi^d,\ldots\rangle \qquad\mbox{where
}\xi^d = \langle\xi _0,\ldots,\xi_d,t+c -[\xi]_d\rangle\in\Gamma_d.
\]
In verifying the last inclusion, note $\sum_{i=0}^d i\dot\xi_i(t) +
(d+1)(1-[\dot\xi(t)]_d)=\break\sum_{i=0}^d (1-[\dot\xi(t)]_i) \leq\sum
_{i\geq0} (1-[\dot\xi(t)]_i) \leq1$. We now argue that $F$ is a
bi-continuous bijection.

Indeed, we first note that $F^{-1}\dvtx\Gamma^*\rightarrow\Gamma^\infty$
is given by
\[
F^{-1}(\varphi) = \langle\varphi^0_0,\ldots, \varphi^d_d, \ldots
\rangle.
\]
In checking $F^{-1}(\varphi)\in\Gamma^\infty$, note for $\varphi
\in
\Gamma^*$ that $\lim_d\sum_{i=0}^d (1-\sum_{l=0}^i\dot\varphi
^l_l(t)) =
\lim_d\sum_{i=0}^d (1-[\dot\varphi^d(t)]_i)\leq1$. Then, by bounded
convergence with respect to the last term in the previous series, $\lim
_d (t +\sum_{i=0}^d c_i -\sum_{i=0}^d \varphi^i_i(t)) = \lim
_d(t+\sum
_{i=0}^d c_i -[\varphi^d(t)]_d) =0$, and so $\sum_{i\geq0} \varphi
^i_i(t) = t+c$. Finally, it is not difficult to see that $F$ and
$F^{-1}$ are both continuous in the product topology.

Now, $\mathcal{X}^{n,\infty}\in\Gamma^*$, $\mathbf X^{n,\infty}\in
\Gamma^\infty$, and $F(\mathbf X^{n,\infty})=\mathcal{X}^{n,\infty}$
for $n\geq1$. Hence,
through the action of $F$, the LDP for $\mathcal{X}^{n,\infty}$
translates to the LDP for $\mathbf X^{n,\infty}$. We
now identify the rate function. Given Propositions
\ref{projectivelimit} and~\ref{rateinfinbounded}, for a degree
distribution $\xi\in\Gamma^\infty$, we identify its rate as
$I^\infty(\xi) = J^\infty(F(\xi))$. Since $\Gamma^\infty$ is closed,
and therefore distributions $\xi\notin\Gamma^\infty$ can never be
attained by $\mathbf{X}^{n,\infty}$, we set $I^\infty(\xi)=\infty$ in
this case. Last, by properties of $F$, as $J^\infty$ is a convex,
good rate function, one obtains readily $I^\infty$ is also a convex,
good rate function.
\end{pf*}

\section{\texorpdfstring{Proof of Corollary \protect\ref{LLN}}{Proof of Corollary 1.7}}\label{sectionLLN}

We verify some properties of $\zeta^d$ in the next lemmas and conclude
the proof of Corollary~\ref{LLN} at the end of the section.
%
%
\begin{lemma}\label{thmzero-cost}
The ODE (\ref{ODEforLLN}) has a unique Carath\'eodory solution
$\zeta^d$.
\end{lemma}
\begin{pf} Any Carath\'eodory solution to ODE (\ref{ODEforLLN}),
given the assumption $p,\beta$ are piecewise continuous, is piecewise
continuously differentiable. Since the defining ODEs are linear, one
can solve them, and so the solution is unique and given by
$\zeta^d = \langle\zeta_0(t),\zeta_1(t),\ldots, \bar\zeta
_{d+1}(t)\rangle$ where,
for $t\in[0,1]$,
%
%
\begin{eqnarray}
\label{zero-cost}\qquad
\zeta_0(t)
&:=& c_0M_0(0,t)+\int_0^t \bigl(1-p(s)\bigr)M_0(s,t) \,ds,\nonumber\\
\zeta_1(t)
&:=& c_1M_1(0,t)\nonumber\\
&&{} +\int_0^t \biggl(p(s) +\bigl(1-p(s)\bigr)\frac{\beta
(s)\zeta
_0(s)}{(1+\beta(s))s+\tilde c+c\beta(s)}\biggr)M_1(s,t) \,ds,
\\
\zeta_i(t)
&:=& c_iM_i(0,t)\nonumber\\
&&{} +\int_0^t \bigl(1-p(s)\bigr)\frac{(i-1+\beta(s))\zeta
_{i-1}(s)}{(1+\beta(s))s+\tilde c+c\beta(s)}M_i(s,t) \,ds\nonumber
\end{eqnarray}
for $2\leq i\leq d$ and
\[
\bar\zeta_{d+1}(t)
:= t+c-\sum_{i=0}^d \zeta_i(t)
= \bar c^d+\int_0^t\bigl(1-p(s)\bigr)\frac{(d+\beta(s))\zeta
_{d}(s)}{(1+\beta
(s))s+\tilde c+c\beta(s)} \,ds.
\]
Here, for $0\leq i\leq d$,
\[
M_i(s,t):= \exp\biggl[-\int_s^t \bigl(1-p(u)\bigr)\frac{i+\beta
(u)}{(1+\beta
(u))u +\tilde c+c\beta(u)} \,du\biggr].
\]
\upqed\end{pf}
%
%
\begin{lemma}\label{lemzero-cost,infinite}
We have $\zeta^d\in\Gamma_d$, and moreover
\[
\sum_{i=0}^\infty\zeta_i(t) = t+c\quad \mbox{and}\quad \sum
_{i=0}^\infty i \zeta_i(t) = t +\tilde c.
\]
\end{lemma}
\begin{pf}
First, from properties of the ODE system and the piecewise continuity
assumption on $p,\beta$ in (ND), $\zeta_i\geq0$, $\zeta_i$ is
Lipschitz with constant $1$ and moreover piecewise continuously
differentiable, and $0\leq[\dot\zeta(t)]_i \leq1$ for $i\geq0$, and
$\sum_{i=0}^d \zeta_i(t) +
\bar\zeta_{d+1}(t)=t+c$ for $d\geq0$
and almost all $t$.
We postpone proving $\sum_{i=0}^d (1-[\dot\zeta(t)]_i) \leq1$ for
$d\geq0$ and a.a. $t$, which would complete the argument to show
$\zeta
^d \in\Gamma_d$, until the end.

We now show $\sum_{i\geq0}\zeta_i(t)=t+c$.
From the defining ODEs (\ref{ODEforLLN}), for $N\geq1$, we have
$1-\sum_{i=0}^N \dot\zeta_i(t) = (1-p(t))\frac{(N+\beta(t))\zeta
_N(t)}{(1+\beta(t))t+\tilde c+c\beta(t)}$,
and hence
%
%
\begin{equation}\label{IEinappendix}\quad
t+\sum_{i=0}^N \zeta_i(0)-\sum_{i=0}^N \zeta_i(t) = \int_0^t
\bigl(1-p(s)\bigr)\frac{(N+\beta(s))\zeta_N(s)}{(1+\beta(s))s+\tilde c+c\beta
(s)} \,ds.
\end{equation}
We obtain, as the integrand on the right-hand side is nonnegative, that
$\sum_{i=0}^N \zeta_i(t)\leq t+\sum_{i=0}^N c_i\leq t+c$ for all $t\geq
0$ and $N\geq1$ where we recall from (LIM) $c=\sum_{i=0}^\infty c_i$.
In particular, $\sum_{i\geq0} \int_0^t \frac{\zeta_i(s)}{s+c}\,ds \leq
t$. Also, the right-hand side of (\ref{IEinappendix}), after a
calculation, is bounded above by $\frac{N+1}{\min\{\beta_0,1\}}\int_0^t
\frac{\zeta_N(s)}{s+c}\,ds$. Hence,\vspace*{1pt} since by nonnegativity
and (LIM) the right-side of (\ref{IEinappendix}) has a limit, this
limit must vanish and $\sum_{i\geq0}\zeta_i(t) = t+c$.

Next, to establish $\sum_{i\geq0}i\zeta_i(t)=t+\tilde c$, again from
the ODEs, for $N\geq1$,
%
%
\begin{eqnarray}\label{ODEweightedsuminappendix}
\sum_{i=0}^N i\dot\zeta_i(t) &=& p(t) + \bigl(1-p(t)\bigr)\frac{\sum
_{i=0}^N(i+\beta(t))\zeta_i(t)}{(1+\beta(t))t+\tilde c+c\beta
(t)}\nonumber\\[-8pt]\\[-8pt]
&&{} - \bigl(1-p(t)\bigr)\frac{(N+1)(N+\beta(t))\zeta_N(t)}{(1+\beta
(t))t+\tilde c+c\beta(t)}.\nonumber
\end{eqnarray}
From\vspace*{1pt} nonnegativity of $\zeta_i$ and $\sum_{i=0}^\infty\zeta_i=t+c$, we
bound the right-hand side of (\ref{ODEweightedsuminappendix}) by
$p(t) + (1-p(t))\frac{\sum_{i=0}^N i\zeta_i(t) + \beta
(t)(t+c)}{(1+\beta(t))t+\tilde c+c\beta(t)}$.
Let\vspace*{1pt} $s_{N}(t):=\sum_{i=0}^N i\zeta_i(t)$. Then,
$\dot s_{N}(t)\leq p(t) + (1-p(t))\frac{s_{N}(t) +
\beta(t)(t+c)}{(1+\beta(t))t+\tilde c+c\beta(t)}$.
Since,
$s_N(t)$ is piecewise continuously differentiable,
we have, by Lemma~\ref{comparingwithsNt}, that $s_N(t) \leq t+
\tilde c$ for $t\geq0$ and $N\geq1$.
Hence,
$\sum_{i=0}^\infty\int_0^t\frac{i\zeta_i(s)}{s+c} \,ds \leq At$ since
$\tilde c\leq Ac$ for some $A>0$ where $\tilde c = \sum_{i\geq
0}ic_i<\infty$.

Now, integrating both sides of ODE (\ref{ODEweightedsuminappendix}), we have
%
%
\begin{eqnarray}\label{IEweightedsuminappendix}
&&
\sum_{i=0}^N i\zeta_i(t) - \sum_{i=0}^N ic_i\nonumber\\
&&\qquad = \int_0^t p(s) \,ds + \int_0^t
\bigl(1-p(s)\bigr)\frac{\sum_{i=0}^N(i+\beta(s))\zeta_i(s)}{(1+\beta
(s))s+\tilde
c+c\beta(s)} \,ds\\
&&\qquad\quad{} - \int_0^t \bigl(1-p(s)\bigr)\frac{(N+1)(N+\beta(s))\zeta
_N(s)}{(1+\beta(s))s+\tilde c+c\beta(s)} \,ds.\nonumber
\end{eqnarray}
From nonnegativity, our estimates and (LIM), the last integral
above has a limit. This last integral in (\ref
{IEweightedsuminappendix}) is bounded above by $\frac{(N+1)^2}{N\min\{
\beta_0,1\}
}\int
_0^t \frac{N\zeta_N(s)}{s+c}\,ds$, and hence its limit must vanish.
Then, using $\sum_{i=0}^\infty\zeta_i(t)= t+ c$, we see $s(t) =
\sum_{i\geq0}i\zeta_i(t)$ satisfies the ODE in Lemma
\ref{comparingwithsNt}, and therefore $s(t) = t + \tilde c$.

Finally, to finish the postponed verification, noting (\ref
{IEinappendix}), we have
\begin{eqnarray*}
\sum_{i=0}^d \bigl(1-[\dot\zeta(t)]_i\bigr) & =& \bigl(1-p(t)\bigr)\frac{s_d(t) + \beta
(t)\sum_{i=0}^d\zeta_i}{(1 + \beta(t))t + \tilde c + c\beta(t)} \\
& \leq&
\frac{ t+\tilde c + \beta(t)(t+c)}{(1 + \beta(t))t + \tilde c +
c\beta
(t)} = 1.
\end{eqnarray*}
\upqed\end{pf}
%
%
\begin{lemma}\label{comparingwithsNt}
The ODE
\[
\dot f(t) = G(t,f(t)) \qquad\mbox{with }G(t,x) = p(t) +
\bigl(1-p(t)\bigr)\frac{x + \beta(t)(t+c)}{(1+\beta(t))t+\tilde c+c\beta(t)}
\]
and initial condition $f(0)=\tilde c$ has unique Carath\'eodory
solution $t+\tilde c$ for $t\geq0$.

In addition, if $u(t)$ is piecewise continuously differentiable,
$u(0)= u_0\leq\tilde c$,
and $\dot u(t) \leq G(t,u(t))$, then $u(t)\leq t+\tilde c$ for $t\geq0$.
\end{lemma}
\begin{pf}
Since the ODE is linear and, from the piecewise continuity assumption
on $p,\beta$ in (ND), $f$ is piecewise
continuously differentiable, we can solve uniquely
\[
f(t) = \tilde c \exp\{B(0,t)\} + \int_0^t \biggl[p(s)+
\frac{(1-p(s))\beta(s)(s+c)}{(1+\beta(s))s + \tilde c +
c\beta(s)}\biggr]
\exp\{B(s,t)\} \,ds,
\]
where $B(q,r) = \int^r_q\frac{1-p(v)}{(1+\beta(v))v + \tilde c +
c\beta
(v)} \,dv$. Recall the convention $0\cdot\infty=0$, so when $c=0$ the
first term
$\tilde c e^{B(0,t)} = 0$ vanishes.
However, $t+\tilde c$ is a solution, and therefore $f(t)$ may be
identified as desired.

The second statement is obtained similarly.
\end{pf}
\begin{pf*}{Proof of Corollary~\ref{LLN}}
Any root of $I_d$ must be a Carath\'eodory solution to ODE (\ref
{ODEforLLN}). Hence, by Lemmas~\ref{thmzero-cost} and \ref
{lemzero-cost,infinite}, $\zeta^d\in\Gamma_d$ is the unique minimizer
of $I_d$. The
LLNs now follow from the LDP upper bound in Theorem~\ref{LDPmain} and
Borel--Cantelli lemma. Statements about ``mass'' and ``weight'' of
$\zeta^\infty$ are proved in Lemma~\ref{lemzero-cost,infinite}.
\end{pf*}

\section{\texorpdfstring{Proof of Corollary \protect\ref{powerlaw}}{Proof of Corollary 1.9}}\label{sectionpowerlaw}

Since $[\zeta^\infty]_i = [\zeta^d]_i$ for $i\leq d$, the proof follows
from the next lemma.
Define, for $o_1, o_2, o_3, o_4, o_5\geq0$, the ODEs,
$O(o_1,o_2,o_3,o_4$, $o_5)$: with initial condition $\varphi(0)=\mathbf c^d$
\begin{eqnarray*}
\dot{\varphi}_0(t)
& = &1-o_1-(1-o_2)\frac{o_3}{1+o_4}\cdot\frac{\varphi
_0(t)}{t+o_5},\\
{}[\dot\varphi(t)]_i
& = &1-(1-o_2)\frac{i+o_3}{1+o_4}\cdot\frac{\varphi_i(t)}{t+o_5}
\qquad\mbox{for } 1\leq i\leq d.
\end{eqnarray*}

One can check that
$\chi(t)$ is the solution to $O(o_1,o_2,o_3,o_4,o_5)$ above for $0\leq
o_2\leq1$, where
%
%
\begin{equation}
\label{soltogeneralODE}\chi_i(t) = b_i (t+o_5) + \sum_{\ell=0}^{i}
a_{i,\ell}\biggl(\frac{o_5}{t+o_5}\biggr)^{(1-o_2)({\ell
+o_3})/({1+o_4})}\qquad
\mbox{for } 0\leq i\leq d.\hspace*{-35pt}
\end{equation}
Here, the sequence $b_i=b_i(o_1,o_2,o_3,o_4,o_5)$ is defined by
$b_0=\frac{1-o_1}{1+(1-o_2){o_3}/({1+o_4})}$,
$b_1=\frac{o_1+(1-o_2){o_3}b_0/({1+o_4})}{1+(1-o_2)
({1+o_3})/({1+o_4})}$,
and, for $i\geq2$,
\begin{eqnarray*}
b_i&=&b_1 \prod_{\ell=2}^{i}\frac{(1-o_2)({\ell
-1+o_3})/({1+o_4})}{1+(1-o_2)({\ell+o_3})/({1+o_4})}\\
&=&b_1\frac{\Gamma(2+o_3+({1+o_4})/({1-o_2}))}{\Gamma(1+o_3)}\frac
{\Gamma(i+o_3)}{\Gamma(i+1+o_3+({1+o_4})/({1-o_2}))}\\
&\sim&\frac{1}{i^{1+({1+o_4})/({1-o_2})}}.
\end{eqnarray*}
The sequence $a_{i,\ell}=a_{i,\ell}(o_1,o_2,o_3,o_4,o_5)$ is given by
$a_{0,0}= c_0-b_0o_5$, and,
for $i\geq1$,
\[
a_{i,\ell}=\frac{i-1+o_3}{i-\ell} a_{i-1,\ell}
\qquad\mbox{where }0\leq\ell<i
\]
and
\[
a_{i,i} = c_i - b_i o_5 - \sum_{\ell=0}^{i-1} a_{i,\ell}.\vadjust{\goodbreak}
\]

Recall now the assumption in Corollary~\ref{powerlaw}:
$0\leq p_{\min}\leq p(\cdot)\leq p_{\max}<1$ and $0<\beta_{\min
}\leq
\beta(\cdot)\leq\beta_{\max}<\infty$.
%
%
\begin{lemma}\label{comparisonlemma} The
systems $O(p_{\min},p_{\max},\beta_{\min},\beta_{\max},\max\{
\tilde c,
c\})$ and $O(p_{\max}$, $p_{\min},\beta_{\max
},\beta_{\min},\min\{\tilde c, c\})$ have
respective unique solutions $\widetilde\zeta$ and $\widehat\zeta$.
Then, for $0\leq i\leq d$ and $t\in[0,1]$, with respect to the
zero-cost trajectory $\zeta^d(t)$ in Corollary~\ref{LLN} with initial
condition $\zeta^d(0)=\mathbf c^d$,
we have
%
\[
[\widehat\zeta(t)]_i \leq[\zeta^d(t)]_i \leq[\widetilde
\zeta(t)]_i.
\]
\end{lemma}
\begin{pf}
The proof that $\widetilde\zeta$ and $\widehat\zeta$ are the unique
solutions uses a similar argument to that in the proof of Lemma \ref
{thmzero-cost}.
We now establish the inequality in the display with respect to
$\widetilde\zeta$ as an analogous proof works for $\widehat\zeta$.
We use induction to see that $[\widetilde\zeta]_i \geq[\zeta]_i$ for
$0\leq i\leq d$.

Since ${\widetilde\zeta}(0)=\zeta(0)=\mathbf c^d$,
from ODEs, $O(p_{\min},p_{\max},\beta_{\min},\beta_{\max
},\max\{\tilde c, c\})$ and (\ref{ODEforLLN}), we have
%
%
\begin{eqnarray}
\label{widezeta0}
\dot{\widetilde\zeta}_0(t) - \dot\zeta_0(t)
& \geq &p(t)-p_{\min}+ (1-p_{\max})\frac{\beta_{\min}(\zeta
_0(t)-\widetilde\zeta_0(t))}{(1+\beta_{\max})(t+\max\{\tilde c, c\})},
\\
\label{widezetai}
\qquad{}[\dot{\widetilde\zeta}(t)]_i - [\dot\zeta(t)]_i
& \geq &(1-p_{\max})\frac{(i+\beta_{\min})(\zeta
_i(t)-\widetilde\zeta
_i(t))}{(1+\beta_{\max})(t+\max\{\tilde c, c\})}.
\end{eqnarray}
For $i=0$, suppose $\widetilde\zeta_0(t)<\zeta_0(t)$ for some $t$.
Then, by continuity, we may assume that
$\widetilde\zeta_0(t)<\zeta_0(t)$ for all $t\in(t_0,t_1]$ for some
$0\leq t_0<t_1\leq1$, and $\widetilde\zeta_0(t_0)=\zeta_0(t_0)$. We
may further arrange $t_0,t_1$, from the piecewise continuity
assumptions in (ND), that $p,\beta$ are continuous on $(t_0,t_1)$.
From the mean value theorem, we find a $t'\in(t_0,t_1)$ such that
$\dot
{\widetilde\zeta}_0(t')<\dot\zeta_0(t')$, which contradicts the ODE
(\ref{widezeta0}) as it gives $\dot{\widetilde\zeta}_0(t') - \dot
\zeta
_0(t')>0$.
Therefore, $\widetilde\zeta_0\geq\zeta_0$.

Now, for $1\leq i\leq d$, suppose $[\widetilde\zeta(t)]_i<[\zeta(t)]_i$
for some $t$.
By induction hypothesis ($[\widetilde\zeta(\cdot)]_{i-1}\geq[\zeta
(\cdot
)]_{i-1}$), we must have $\widetilde\zeta_i(t)<\zeta_i(t)$.
Since $[\widetilde\zeta(\cdot)]_i$, $[\zeta(\cdot)]_i$,
$\widetilde\zeta
_i(\cdot)$ and $\zeta_i(\cdot)$ are continuous functions, as for the
case $i=0$,
we may assume $[\widetilde\zeta(t)]_i< [\zeta(t)]_i$ and
$\widetilde\zeta_i(t)<\zeta_i(t)$, and $p,\beta$ are continuous for all
$t\in(t_0,t_1)$ for some $0\leq t_0<t_1\leq1$, and also $\widetilde
\zeta_i(t_0)=\zeta_i(t_0)$.
By the mean value theorem for $[\widetilde\zeta(t)]_i-[\zeta(t)]_i$,
there is $t'\in(t_0,t_1)$ such that $[\dot{\widetilde\zeta
}(t')]_i<[\dot
\zeta(t')]_i$.
But (\ref{widezetai}) gives $[\dot{\widetilde\zeta}(t')]_i -
[\dot
\zeta(t')]_i>0$,
a contradiction. Therefore $[\widetilde\zeta]_i\geq[\zeta]_i$.
\end{pf}
\begin{pf*}{Proof of Corollary~\ref{powerlaw}}
Given Lemma~\ref{comparisonlemma}, we need only detail the solutions
$\widetilde\zeta$ and $\widehat\zeta$ when the initial
configuration is
``small'' and ``large,'' respectively. To this end, when the initial
configuration is ``small'' ($c_i\equiv0$), $\widetilde\zeta$,
$\widehat\zeta$ are linear, namely $\widetilde\zeta_i(t)=\widetilde b_i
t$, and $\widehat\zeta_i(t)=\widehat b_i t$, where $\widetilde
b_i:=b_i(p_{\min}, p_{\max}, \beta_{\min},\break \beta_{\max}, 0)$ and
$\widehat b_i:= b_i(p_{\max}, p_{\min},
\beta_{\max}, \beta_{\min}, 0)$ [cf. (\ref{soltogeneralODE})].

On the other hand, when the initial configuration is ``large'' ($c_i>0$
for some $0\leq i\leq d+1$), as $t\uparrow\infty$,
$\widetilde\zeta_i(t)=(\widetilde b_i + o(1)) t$ and $\widehat\zeta
_i(t)=(\widehat b_i + o(1)) t$.
\end{pf*}

\section*{Acknowledgments}

We thank the referees/editors for constructive comments which helped
improve the exposition of the paper.


%

\printaddresses


\begin{thebibliography}{53}

\bibitem{Albert-Barabasi-02}
%
\begin{barticle}[mr]
\bauthor{\bsnm{Albert},~\bfnm{R{\'e}ka}\binits{R.}} \AND
\bauthor{\bsnm{Barab{\'a}si},~\bfnm{Albert-L{\'a}szl{\'o}}\binits{A.-L.}}
(\byear{2002}).
\btitle{Statistical mechanics of complex networks}.
\bjournal{Rev. Modern Phys.}
\bvolume{74}
\bpages{47--97}.
\bid{doi={10.1103/RevModPhys.74.47}, issn={0034-6861}, mr={1895096}}
\bptok{imsref}%
\end{barticle}
%
\endbibitem

\bibitem{AGS08}
%
\begin{barticle}[mr]
\bauthor{\bsnm{Athreya},~\bfnm{Krishna~B.}\binits{K.~B.}},
\bauthor{\bsnm{Ghosh},~\bfnm{Arka~P.}\binits{A.~P.}} \AND
\bauthor{\bsnm{Sethuraman},~\bfnm{Sunder}\binits{S.}}
(\byear{2008}).
\btitle{Growth of preferential attachment random graphs via continuous-time
branching processes}.
\bjournal{Proc. Indian Acad. Sci. Math. Sci.}
\bvolume{118}
\bpages{473--494}.
\bid{doi={10.1007/s12044-008-0036-2}, issn={0253-4142}, mr={2450248}}
\bptok{imsref}%
\end{barticle}
%
\endbibitem

\bibitem{Barabasi}
%
\begin{barticle}[mr]
\bauthor{\bsnm{Barab{\'a}si},~\bfnm{Albert-L{\'a}szl{\'o}}\binits{A.-L.}}
(\byear{2009}).
\btitle{Scale-free networks: A decade and beyond}.
\bjournal{Science}
\bvolume{325}
\bpages{412--413}.
\bid{doi={10.1126/science.1173299}, issn={0036-8075}, mr={2548299}}
\bptok{imsref}%
\end{barticle}
%
\endbibitem

\bibitem{Albert-Barabasi-99}
%
\begin{barticle}[mr]
\bauthor{\bsnm{Barab{\'a}si},~\bfnm{Albert-L{\'a}szl{\'o}}\binits
{A.-L.}} \AND
\bauthor{\bsnm{Albert},~\bfnm{R{\'e}ka}\binits{R.}}
(\byear{1999}).
\btitle{Emergence of scaling in random networks}.
\bjournal{Science}
\bvolume{286}
\bpages{509--512}.
\bid{doi={10.1126/science.286.5439.509}, issn={0036-8075}, mr={2091634}}
\bptok{imsref}%
\end{barticle}
%
\endbibitem

\bibitem{BBV}
%
\begin{barticle}[auto:STB|2012/04/30|08:06:40]
\bauthor{\bsnm{Barrat},~\bfnm{A.}\binits{A.}},
\bauthor{\bsnm{Barthelemy},~\bfnm{M.}\binits{M.}},
\bauthor{\bsnm{Pastor-Satorras},~\bfnm{R.}\binits{R.}} \AND
\bauthor{\bsnm{Vespignani},~\bfnm{A.}\binits{A.}}
(\byear{2004}).
\btitle{The architecture of complex weighted networks}.
\bjournal{PNAS}
\bvolume{101}
\bpages{3747--3752}.
\bptok{imsref}%
\end{barticle}
%
\endbibitem

\bibitem{Berger-Borgs-Chayes-Saberi-09}
%
\begin{bmisc}[auto:STB|2012/04/30|08:06:40]
\bauthor{\bsnm{Berger},~\bfnm{N.}\binits{N.}},
\bauthor{\bsnm{Borgs},~\bfnm{C.}\binits{C.}},
\bauthor{\bsnm{Chayes},~\bfnm{J.}\binits{J.}} \AND
\bauthor{\bsnm{Saberi},~\bfnm{A.}\binits{A.}}
(\byear{2009}).
\bhowpublished{A weak local limit for preferential attachment graphs.
Preprint}.
\bptok{imsref}%
\end{bmisc}
%
\endbibitem

\bibitem{BBCS-05}
%
\begin{binproceedings}[mr]
\bauthor{\bsnm{Berger},~\bfnm{Noam}\binits{N.}},
\bauthor{\bsnm{Borgs},~\bfnm{Christian}\binits{C.}},
\bauthor{\bsnm{Chayes},~\bfnm{Jennifer~T.}\binits{J.~T.}} \AND
\bauthor{\bsnm{Saberi},~\bfnm{Amin}\binits{A.}}
(\byear{2005}).
\btitle{On the spread of viruses on the internet}.
In \bbooktitle{Proceedings of the {S}ixteenth {A}nnual {ACM}--{SIAM} {S}ymposium
on {D}iscrete {A}lgorithms}
\bpages{301--310 (electronic)}.
\bpublisher{ACM}, \baddress{New York}.
\bid{mr={2298278}}
\bptok{imsref}%
\end{binproceedings}
%
\endbibitem

\bibitem{Bhamidi}
%
\begin{bmisc}[auto:STB|2012/04/30|08:06:40]
\bauthor{\bsnm{Bhamidi},~\bfnm{S.}\binits{S.}}
(\byear{2007}).
\bhowpublished{Universal techniques to analyze preferential attachment trees:
Global and local analysis. Preprint. Available at
\texttt{%
\href{http://www.unc.edu/\textasciitilde bhamidi/preferent.pdf}{http://www.unc.edu/}
\href{http://www.unc.edu/\textasciitilde bhamidi/preferent.pdf}{\textasciitilde bhamidi/preferent.pdf}}.}
\bptok{imsref}%
\end{bmisc}
%
\endbibitem

\bibitem{Bollobas-Riordan}
%
\begin{barticle}[mr]
\bauthor{\bsnm{Bollob{\'a}s},~\bfnm{B{\'e}la}\binits{B.}} \AND
\bauthor{\bsnm{Riordan},~\bfnm{Oliver}\binits{O.}}
(\byear{2004}).
\btitle{The diameter of a scale-free random graph}.
\bjournal{Combinatorica}
\bvolume{24}
\bpages{5--34}.
\bid{doi={10.1007/s00493-004-0002-2}, issn={0209-9683}, mr={2057681}}
\bptok{imsref}%
\end{barticle}
%
\endbibitem

\bibitem{BRST01}
%
\begin{barticle}[mr]
\bauthor{\bsnm{Bollob{\'a}s},~\bfnm{B{\'e}la}\binits{B.}},
\bauthor{\bsnm{Riordan},~\bfnm{Oliver}\binits{O.}},
\bauthor{\bsnm{Spencer},~\bfnm{Joel}\binits{J.}} \AND
\bauthor{\bsnm{Tusn{\'a}dy},~\bfnm{G{\'a}bor}\binits{G.}}
(\byear{2001}).
\btitle{The degree sequence of a scale-free random graph process}.
\bjournal{Random Structures Algorithms}
\bvolume{18}
\bpages{279--290}.
\bid{doi={10.1002/rsa.1009}, issn={1042-9832}, mr={1824277}}
\bptok{imsref}%
\end{barticle}
%
\endbibitem

\bibitem{BCLSV}
%
\begin{barticle}[mr]
\bauthor{\bsnm{Borgs},~\bfnm{Christian}\binits{C.}},
\bauthor{\bsnm{Chayes},~\bfnm{Jennifer}\binits{J.}},
\bauthor{\bsnm{Lov{\'a}sz},~\bfnm{L{\'a}szl{\'o}}\binits{L.}},
\bauthor{\bsnm{S{\'o}s},~\bfnm{Vera}\binits{V.}} \AND
\bauthor{\bsnm{Vesztergombi},~\bfnm{Katalin}\binits{K.}}
(\byear{2011}).
\btitle{Limits of randomly grown graph sequences}.
\bjournal{European J. Combin.}
\bvolume{32}
\bpages{985--999}.
\bid{doi={10.1016/j.ejc.2011.03.015}, issn={0195-6698}, mr={2825531}}
\bptok{imsref}%
\end{barticle}
%
\endbibitem

\bibitem{BS}
%
\begin{bbook}[mr]
\beditor{\bsnm{Bornholdt},~\bfnm{Stefan}\binits{S.}} \AND
\beditor{\bsnm{Schuster},~\bfnm{Heinz~Georg}\binits{H.~G.}}, eds.
(\byear{2003}).
\btitle{Handbook of Graphs and Networks: From the Genome to the Internet}.
\bpublisher{Wiley--VCH}, \baddress{Weinheim}.
\bid{mr={2016116}}
\bptok{imsref}%
\end{bbook}
%
\endbibitem

\bibitem{Bryc-Minda-S}
%
\begin{barticle}[mr]
\bauthor{\bsnm{Bryc},~\bfnm{Wlodek}\binits{W.}},
\bauthor{\bsnm{Minda},~\bfnm{David}\binits{D.}} \AND
\bauthor{\bsnm{Sethuraman},~\bfnm{Sunder}\binits{S.}}
(\byear{2009}).
\btitle{Large deviations for the leaves in some random trees}.
\bjournal{Adv. in Appl. Probab.}
\bvolume{41}
\bpages{845--873}.
\bid{doi={10.1239/aap/1253281066}, issn={0001-8678}, mr={2571319}}
\bptok{imsref}%
\end{barticle}
%
\endbibitem

\bibitem{Cald}
%
\begin{bbook}[auto:STB|2012/04/30|08:06:40]
\bauthor{\bsnm{Caldarelli},~\bfnm{Guido}\binits{G.}}
(\byear{2007}).
\btitle{Scale-Free Networks: Complex Webs in Nature and Technology}.
\bpublisher{Oxford Univ. Press}, \baddress{Oxford}.
\bptok{imsref}%
\end{bbook}
%
\endbibitem

\bibitem{Chung-Handjani-Jungreis}
%
\begin{barticle}[mr]
\bauthor{\bsnm{Chung},~\bfnm{Fan}\binits{F.}},
\bauthor{\bsnm{Handjani},~\bfnm{Shirin}\binits{S.}} \AND
\bauthor{\bsnm{Jungreis},~\bfnm{Doug}\binits{D.}}
(\byear{2003}).
\btitle{Generalizations of {P}olya's urn problem}.
\bjournal{Ann. Comb.}
\bvolume{7}
\bpages{141--153}.
\bid{doi={10.1007/s00026-003-0178-y}, issn={0218-0006}, mr={1994572}}
\bptok{imsref}%
\end{barticle}
%
\endbibitem

\bibitem{CL06}
%
\begin{bbook}[mr]
\bauthor{\bsnm{Chung},~\bfnm{Fan}\binits{F.}} \AND
\bauthor{\bsnm{Lu},~\bfnm{Linyuan}\binits{L.}}
(\byear{2006}).
\btitle{Complex Graphs and Networks}.
\bseries{CBMS Regional Conference Series in Mathematics}
\bvolume{107}.
\bpublisher{Amer. Math. Soc.},
\baddress{Providence, RI}.
\bid{mr={2248695}}
\bptok{imsref}%
\end{bbook}
%
\endbibitem

\bibitem{CH}
%
\begin{bbook}[auto:STB|2012/04/30|08:06:40]
\bauthor{\bsnm{Cohen},~\bfnm{R.}\binits{R.}} \AND
\bauthor{\bsnm{Havlin},~\bfnm{S.}\binits{S.}}
(\byear{2010}).
\btitle{Complex Networks: Structure Robustness and Function}.
\bpublisher{Cambridge Univ. Press}, \baddress{Cambridge}.
\bptok{imsref}%
\end{bbook}
%
\endbibitem

\bibitem{Cooper-Frieze-general}
%
\begin{barticle}[mr]
\bauthor{\bsnm{Cooper},~\bfnm{Colin}\binits{C.}} \AND
\bauthor{\bsnm{Frieze},~\bfnm{Alan}\binits{A.}}
(\byear{2003}).
\btitle{A general model of web graphs}.
\bjournal{Random Structures Algorithms}
\bvolume{22}
\bpages{311--335}.
\bid{doi={10.1002/rsa.10084}, issn={1042-9832}, mr={1966545}}
\bptok{imsref}%
\end{barticle}
%
\endbibitem

\bibitem{Cooper-Frieze07}
%
\begin{barticle}[mr]
\bauthor{\bsnm{Cooper},~\bfnm{Colin}\binits{C.}} \AND
\bauthor{\bsnm{Frieze},~\bfnm{Alan}\binits{A.}}
(\byear{2007}).
\btitle{The cover time of the preferential attachment graph}.
\bjournal{J.~Combin. Theory Ser. B}
\bvolume{97}
\bpages{269--290}.
\bid{doi={10.1016/j.jctb.2006.05.007}, issn={0095-8956}, mr={2290325}}
\bptok{imsref}%
\end{barticle}
%
\endbibitem

\bibitem{Dembo-Zeitouni}
%
\begin{bbook}[mr]
\bauthor{\bsnm{Dembo},~\bfnm{Amir}\binits{A.}} \AND
\bauthor{\bsnm{Zeitouni},~\bfnm{Ofer}\binits{O.}}
(\byear{1998}).
\btitle{Large Deviations Techniques and Applications},
\bedition{2nd} ed.
\bseries{Applications of Mathematics (New York)}
\bvolume{38}.
\bpublisher{Springer}, \baddress{New York}.
\bid{mr={1619036}}
\bptok{imsref}%
\end{bbook}
%
\endbibitem

\bibitem{dereich-2009}
%
\begin{barticle}[mr]
\bauthor{\bsnm{Dereich},~\bfnm{Steffen}\binits{S.}} \AND
\bauthor{\bsnm{M{\"o}rters},~\bfnm{Peter}\binits{P.}}
(\byear{2009}).
\btitle{Random networks with sublinear preferential attachment: Degree
evolutions}.
\bjournal{Electron. J. Probab.}
\bvolume{14}
\bpages{1222--1267}.
\bid{doi={10.1214/EJP.v14-647}, issn={1083-6489}, mr={2511283}}
\bptok{imsref}%
\end{barticle}
%
\endbibitem

\bibitem{Dereich-Monch-Morters}
%
\begin{bmisc}[mr]
\bauthor{\bsnm{Dereich},~\bfnm{Steffen}\binits{S.}},
\bauthor{\bsnm{M\"onch},~\bfnm{C.}\binits{C.}} \AND
\bauthor{\bsnm{M{\"o}rters},~\bfnm{Peter}\binits{P.}}
(\byear{2011}).
\bhowpublished{Typical distances in ultrasmall random networks.
Available at arXiv:\arxivurl{1102.5680v1}.}
\bptok{imsref}%
\end{bmisc}
%
\endbibitem

\bibitem{DKM}
%
\begin{barticle}[mr]
\bauthor{\bsnm{Dorogovtsev},~\bfnm{S.~N.}\binits{S.~N.}},
\bauthor{\bsnm{Krapivsky},~\bfnm{P.~L.}\binits{P.~L.}} \AND
\bauthor{\bsnm{Mend{\'e}s},~\bfnm{J.~F.~F.}\binits{J.~F.~F.}}
(\byear{2008}).
\btitle{Transition from small to large world in growing networks}.
\bjournal{Europhys. Lett. EPL}
\bvolume{81}
\bpages{Art. 30004, 5}.
\bid{doi={10.1209/0295-5075/81/30004}, issn={0295-5075}, mr={2443959}}
\bptok{imsref}%
\end{barticle}
%
\endbibitem

\bibitem{DMpaper}
%
\begin{barticle}[auto:STB|2012/04/30|08:06:40]
\bauthor{\bsnm{Dorogovtsev},~\bfnm{S.~N.}\binits{S.~N.}} \AND
\bauthor{\bsnm{Mendes},~\bfnm{J.~F.~F.}\binits{J.~F.~F.}}
(\byear{2000}).
\btitle{Evolution of networks with aging of sites}.
\bjournal{Phys. Rev. E}
\bvolume{62}
\bpages{1842--1845}.
\bptok{imsref}%
\end{barticle}
%
\endbibitem

\bibitem{Dorogovtsev-Mendes}
%
\begin{barticle}[auto:STB|2012/04/30|08:06:40]
\bauthor{\bsnm{Dorogovtsev},~\bfnm{S.~N.}\binits{S.~N.}} \AND
\bauthor{\bsnm{Mendes},~\bfnm{J.~F.~F.}\binits{J.~F.~F.}}
(\byear{2001}).
\btitle{Scaling properties of scale-free evolving networks: Continuous
approach}.
\bjournal{Phys. Rev. E}
\bvolume{63}
\bpages{056125 19 pp.}
\bptok{imsref}%
\end{barticle}
%
\endbibitem

\bibitem{DM}
%
\begin{bbook}[mr]
\bauthor{\bsnm{Dorogovtsev},~\bfnm{S.~N.}\binits{S.~N.}} \AND
\bauthor{\bsnm{Mendes},~\bfnm{J.~F.~F.}\binits{J.~F.~F.}}
(\byear{2003}).
\btitle{Evolution of Networks: From Biological Nets to the Internet and WWW}.
\bpublisher{Oxford Univ. Press}, \baddress{Oxford}.
\bid{doi={10.1093/acprof:oso/9780198515906.001.0001}, mr={1993912}}
\bptok{imsref}%
\end{bbook}
%
\endbibitem

\bibitem{DEM}
%
\begin{bmisc}[auto:STB|2012/04/30|08:06:40]
\bauthor{\bsnm{Drinea},~\bfnm{E.}\binits{E.}},
\bauthor{\bsnm{Enachescu},~\bfnm{M.}\binits{M.}} \AND
\bauthor{\bsnm{Mitzenmacher},~\bfnm{M.}\binits{M.}}
(\byear{2001}).
\bhowpublished{Variations on random graph models for the web. Harvard Technical
Report TR-06-01.}
\bptok{imsref}%
\end{bmisc}
%
\endbibitem

\bibitem{Drinea-Frieze-Mitzenmacher}
%
\begin{binproceedings}[auto:STB|2012/04/30|08:06:40]
\bauthor{\bsnm{Drinea},~\bfnm{E.}\binits{E.}},
\bauthor{\bsnm{Frieze},~\bfnm{A.}\binits{A.}} \AND
\bauthor{\bsnm{Mitzenmacher},~\bfnm{M.}\binits{M.}}
(\byear{2002}).
\btitle{Balls and Bins models with feedback}.
In \bbooktitle{Proc. of the 11th ACM--SIAM Symposium on Discrete Algorithms
(SODA)}
\bpages{308--315}.
\bpublisher{SIAM},
\baddress{Philadelphia, PA}.
\bptok{imsref}%
\end{binproceedings}
%
\endbibitem

\bibitem{Dupuis-Ellis}
%
\begin{bbook}[mr]
\bauthor{\bsnm{Dupuis},~\bfnm{Paul}\binits{P.}} \AND
\bauthor{\bsnm{Ellis},~\bfnm{Richard~S.}\binits{R.~S.}}
(\byear{1997}).
\btitle{A Weak Convergence Approach to the Theory of Large Deviations}.
\bpublisher{Wiley}, \baddress{New York}.
\bid{doi={10.1002/9781118165904}, mr={1431744}}
\bptok{imsref}%
\end{bbook}
%
\endbibitem

\bibitem{Du06}
%
\begin{bbook}[mr]
\bauthor{\bsnm{Durrett},~\bfnm{Rick}\binits{R.}}
(\byear{2007}).
\btitle{Random Graph Dynamics}.
\bpublisher{Cambridge Univ. Press}, \baddress{Cambridge}.
\bid{mr={2271734}}
\bptnote{check year}%
\bptok{imsref}%
\end{bbook}
%
\endbibitem

\bibitem{Fort}
%
\begin{barticle}[auto:STB|2012/04/30|08:06:40]
\bauthor{\bsnm{Fortunato},~\bfnm{S.}\binits{S.}},
\bauthor{\bsnm{Flammini},~\bfnm{A.}\binits{A.}} \AND
\bauthor{\bsnm{Menczer},~\bfnm{F.}\binits{F.}}
(\byear{2006}).
\btitle{Scale-free network growth by ranking}.
\bjournal{Phys. Rev. Lett.}
\bvolume{96}
\bpages{218701-1--218701-4}.
\bptok{imsref}%
\end{barticle}
%
\endbibitem

\bibitem{FVC}
%
\begin{barticle}[mr]
\bauthor{\bsnm{Frieze},~\bfnm{Alan}\binits{A.}},
\bauthor{\bsnm{Vera},~\bfnm{Juan}\binits{J.}} \AND
\bauthor{\bsnm{Chakrabarti},~\bfnm{Soumen}\binits{S.}}
(\byear{2006}).
\btitle{The influence of search engines on preferential attachment}.
\bjournal{Internet Math.}
\bvolume{3}
\bpages{361--381}.
\bid{issn={1542-7951}, mr={2372548}}
\bptok{imsref}%
\end{barticle}
%
\endbibitem

\bibitem{Gjoka-Kurant-Butts-Markopoulou}
%
\begin{barticle}[auto:STB|2012/04/30|08:06:40]
\bauthor{\bsnm{Gjoka},~\bfnm{M.}\binits{M.}},
\bauthor{\bsnm{Kurant},~\bfnm{M.}\binits{M.}},
\bauthor{\bsnm{Butts},~\bfnm{C.~T.}\binits{C.~T.}} \AND
\bauthor{\bsnm{Markopoulou},~\bfnm{A.}\binits{A.}}
(\byear{2011}).
\btitle{Practical recommendations on crawling online social networks}.
\bjournal{IEEE Journal on Selected Areas in Communications}
\bvolume{29}
\bpages{1872--1892}.
\bptok{imsref}%
\end{barticle}
%
\endbibitem

\bibitem{Jans}
%
\begin{barticle}[mr]
\bauthor{\bsnm{Janssen},~\bfnm{Jeannette}\binits{J.}} \AND
\bauthor{\bsnm{Pra{\l}at},~\bfnm{Pawe{\l}}\binits{P.}}
(\byear{2010}).
\btitle{Rank-based attachment leads to power law graphs}.
\bjournal{SIAM J. Discrete Math.}
\bvolume{24}
\bpages{420--440}.
\bid{doi={10.1137/080716967}, issn={0895-4801}, mr={2646095}}
\bptok{imsref}%
\end{barticle}
%
\endbibitem

\bibitem{Katona}
%
\begin{barticle}[mr]
\bauthor{\bsnm{Katona},~\bfnm{Zsolt}\binits{Z.}}
(\byear{2005}).
\btitle{Width of a scale-free tree}.
\bjournal{J. Appl. Probab.}
\bvolume{42}
\bpages{839--850}.
\bid{issn={0021-9002}, mr={2157524}}
\bptok{imsref}%
\end{barticle}
%
\endbibitem

\bibitem{KR}
%
\begin{barticle}[auto:STB|2012/04/30|08:06:40]
\bauthor{\bsnm{Krapivsky},~\bfnm{P.}\binits{P.}} \AND
\bauthor{\bsnm{Redner},~\bfnm{S.}\binits{S.}}
(\byear{2001}).
\btitle{Organization of growing random networks}.
\bjournal{Phys. Rev. E}
\bvolume{63}
\bpages{066123-1--066123-14}.
\bptok{imsref}%
\end{barticle}
%
\endbibitem

\bibitem{KR-finite}
%
\begin{barticle}[mr]
\bauthor{\bsnm{Krapivsky},~\bfnm{P.~L.}\binits{P.~L.}} \AND
\bauthor{\bsnm{Redner},~\bfnm{S.}\binits{S.}}
(\byear{2002}).
\btitle{Finiteness and fluctuations in growing networks}.
\bjournal{J. Phys. A}
\bvolume{35}
\bpages{9517--9534}.
\bid{doi={10.1088/0305-4470/35/45/302}, issn={0305-4470}, mr={1946936}}
\bptok{imsref}%
\end{barticle}
%
\endbibitem

\bibitem{KRR}
%
\begin{barticle}[pbm]
\bauthor{\bsnm{Krapivsky},~\bfnm{P.~L.}\binits{P.~L.}},
\bauthor{\bsnm{Rodgers},~\bfnm{G.~J.}\binits{G.~J.}} \AND
\bauthor{\bsnm{Redner},~\bfnm{S.}\binits{S.}}
(\byear{2001}).
\btitle{Degree distributions of growing networks}.
\bjournal{Phys. Rev. Lett.}
\bvolume{86}
\bpages{5401--5404}.
\bid{issn={0031-9007}, pmid={11384508}}
\bptok{imsref}%
\end{barticle}
%
\endbibitem

\bibitem{MPS}
%
\begin{barticle}[mr]
\bauthor{\bsnm{Mihail},~\bfnm{Milena}\binits{M.}},
\bauthor{\bsnm{Papadimitriou},~\bfnm{Christos}\binits{C.}} \AND
\bauthor{\bsnm{Saberi},~\bfnm{Amin}\binits{A.}}
(\byear{2006}).
\btitle{On certain connectivity properties of the internet topology}.
\bjournal{J. Comput. System Sci.}
\bvolume{72}
\bpages{239--251}.
\bid{doi={10.1016/j.jcss.2005.06.009}, issn={0022-0000}, mr={2205286}}
\bptok{imsref}%
\end{barticle}
%
\endbibitem

\bibitem{Mitzenmacher}
%
\begin{barticle}[mr]
\bauthor{\bsnm{Mitzenmacher},~\bfnm{Michael}\binits{M.}}
(\byear{2004}).
\btitle{A brief history of generative models for power law and lognormal
distributions}.
\bjournal{Internet Math.}
\bvolume{1}
\bpages{226--251}.
\bid{issn={1542-7951}, mr={2077227}}
\bptok{imsref}%
\end{barticle}
%
\endbibitem

\bibitem{Mori-01}
%
\begin{barticle}[mr]
\bauthor{\bsnm{M{\'o}ri},~\bfnm{T.~F.}\binits{T.~F.}}
(\byear{2002}).
\btitle{On random trees}.
\bjournal{Studia Sci. Math. Hungar.}
\bvolume{39}
\bpages{143--155}.
\bid{doi={10.1556/SScMath.39.2002.1-2.9}, issn={0081-6906}, mr={1909153}}
\bptok{imsref}%
\end{barticle}
%
\endbibitem

\bibitem{Mori-05}
%
\begin{barticle}[mr]
\bauthor{\bsnm{M{\'o}ri},~\bfnm{Tam{\'a}s~F.}\binits{T.~F.}}
(\byear{2005}).
\btitle{The maximum degree of the {B}arab\'asi--{A}lbert random tree}.
\bjournal{Combin. Probab. Comput.}
\bvolume{14}
\bpages{339--348}.
\bid{doi={10.1017/S0963548304006133}, issn={0963-5483}, mr={2138118}}
\bptok{imsref}%
\end{barticle}
%
\endbibitem

\bibitem{NW}
%
\begin{bbook}[mr]
\beditor{\bsnm{Newman},~\bfnm{Mark}\binits{M.}},
\beditor{\bsnm{Barab{\'a}si},~\bfnm{Albert-L{\'a}szl{\'o}}\binits{A.-L.}}
\AND\beditor{\bsnm{Watts},~\bfnm{Duncan~J.}\binits{D.~J.}}, eds.
(\byear{2006}).
\btitle{The Structure and Dynamics of Networks}.
\bpublisher{Princeton Univ. Press}, \baddress{Princeton, NJ}.
\bid{mr={2352222}}
\bptok{imsref}%
\end{bbook}
%
\endbibitem

\bibitem{Newman}
%
\begin{barticle}[mr]
\bauthor{\bsnm{Newman},~\bfnm{M.~E.~J.}\binits{M.~E.~J.}}
(\byear{2003}).
\btitle{The structure and function of complex networks}.
\bjournal{SIAM Rev.}
\bvolume{45}
\bpages{167--256 (electronic)}.
\bid{doi={10.1137/S003614450342480}, issn={0036-1445}, mr={2010377}}
\bptok{imsref}%
\end{barticle}
%
\endbibitem

\bibitem{Newman10}
%
\begin{bbook}[mr]
\bauthor{\bsnm{Newman},~\bfnm{M.~E.~J.}\binits{M.~E.~J.}}
(\byear{2010}).
\btitle{Networks: An Introduction}.
\bpublisher{Oxford Univ. Press}, \baddress{Oxford}.
\bid{doi={10.1093/acprof:oso/9780199206650.001.0001}, mr={2676073}}
\bptok{imsref}%
\end{bbook}
%
\endbibitem

\bibitem{Oliveira-Spencer05}
%
\begin{barticle}[mr]
\bauthor{\bsnm{Oliveira},~\bfnm{Roberto}\binits{R.}} \AND
\bauthor{\bsnm{Spencer},~\bfnm{Joel}\binits{J.}}
(\byear{2005}).
\btitle{Connectivity transitions in networks with super-linear preferential
attachment}.
\bjournal{Internet Math.}
\bvolume{2}
\bpages{121--163}.
\bid{issn={1542-7951}, mr={2193157}}
\bptok{imsref}%
\end{barticle}
%
\endbibitem

\bibitem{RS}
%
\begin{bmisc}[auto:STB|2012/04/30|08:06:40]
\bauthor{\bsnm{R{\'a}th},~\bfnm{B.}\binits{B.}} \AND
\bauthor{\bsnm{Szak{\'a}cs},~\bfnm{L.}\binits{L.}}
(\byear{2011}).
\bhowpublished{Multigraph limit of the dense configuration model and the
preferential attachment graph. Preprint. Available at arXiv:\arxivurl{1106.2058}.}
\bptok{imsref}%
\end{bmisc}
%
\endbibitem

\bibitem{Rudas-Toth-Valko-07}
%
\begin{barticle}[mr]
\bauthor{\bsnm{Rudas},~\bfnm{Anna}\binits{A.}},
\bauthor{\bsnm{T{\'o}th},~\bfnm{B{\'a}lint}\binits{B.}} \AND
\bauthor{\bsnm{Valk{\'o}},~\bfnm{Benedek}\binits{B.}}
(\byear{2007}).
\btitle{Random trees and general branching processes}.
\bjournal{Random Structures Algorithms}
\bvolume{31}
\bpages{186--202}.
\bid{doi={10.1002/rsa.20137}, issn={1042-9832}, mr={2343718}}
\bptok{imsref}%
\end{barticle}
%
\endbibitem

\bibitem{Simkin}
%
\begin{barticle}[mr]
\bauthor{\bsnm{Simkin},~\bfnm{M.~V.}\binits{M.~V.}} \AND
\bauthor{\bsnm{Roychowdhury},~\bfnm{V.~P.}\binits{V.~P.}}
(\byear{2011}).
\btitle{Re-inventing {W}illis}.
\bjournal{Phys. Rep.}
\bvolume{502}
\bpages{1--35}.
\bid{doi={10.1016/j.physrep.2010.12.004}, issn={0370-1573}, mr={2788549}}
\bptok{imsref}%
\end{barticle}
%
\endbibitem

\bibitem{Simon}
%
\begin{barticle}[mr]
\bauthor{\bsnm{Simon},~\bfnm{Herbert~A.}\binits{H.~A.}}
(\byear{1955}).
\btitle{On a class of skew distribution functions}.
\bjournal{Biometrika}
\bvolume{42}
\bpages{425--440}.
\bid{issn={0006-3444}, mr={0073085}}
\bptok{imsref}%
\end{barticle}
%
\endbibitem

\bibitem{Yule}
%
\begin{barticle}[auto:STB|2012/04/30|08:06:40]
\bauthor{\bsnm{Yule},~\bfnm{G.~U.}\binits{G.~U.}}
(\byear{1924}).
\btitle{A mathematical theory of evolution, based on the conclusions of Dr.
J.~C. Willis}.
\bjournal{Philos. Trans. Roy. Soc. London Ser. B}
\bvolume{213}
\bpages{21--87}.
\bptok{imsref}%
\end{barticle}
%
\endbibitem

\bibitem{Zhang-Dupuis-08}
%
\begin{barticle}[mr]
\bauthor{\bsnm{Zhang},~\bfnm{Jim~X.}\binits{J.~X.}} \AND
\bauthor{\bsnm{Dupuis},~\bfnm{Paul}\binits{P.}}
(\byear{2008}).
\btitle{Large-deviation approximations for general occupancy models}.
\bjournal{Combin. Probab. Comput.}
\bvolume{17}
\bpages{437--470}.
\bid{doi={10.1017/S0963548307008681}, issn={0963-5483}, mr={2410397}}
\bptok{imsref}%
\end{barticle}
%
\endbibitem

\end{thebibliography}
\end{document}